\documentclass[lettersize,journal]{IEEEtran}
\usepackage{amsmath,amsfonts}
\usepackage{algorithmic}

\usepackage{cite}
\usepackage{amsmath,amssymb,amsfonts}
\usepackage{algorithmic}
\usepackage{graphicx}
\usepackage{algorithm,algorithmic}
\usepackage{hyperref}
\hypersetup{hidelinks=true}
\usepackage{textcomp}
\usepackage{mathrsfs}  
\usepackage{array}
\usepackage{multirow}
\usepackage{booktabs}

\newcommand{\tabincell}[2]{\begin{tabular}{@{}#1@{}}#2\end{tabular}} 

\usepackage{mathrsfs}

\usepackage{array}
\usepackage[caption=false,font=normalsize,labelfont=sf,textfont=sf]{subfig}
\usepackage{textcomp}
\usepackage{stfloats}
\usepackage{url}
\usepackage{verbatim}
\usepackage{graphicx}
\def\BibTeX{{\rm B\kern-.05em{\sc i\kern-.025em b}\kern-.08em
    T\kern-.1667em\lower.7ex\hbox{E}\kern-.125emX}}
\usepackage{balance}
\begin{document}
\title{A distributed augmented Lagrangian decomposition algorithm for constrained optimization}
\author{
Wenyou Guo, Ting Qu, Hainan Huang, and  Yafeng Wei
\thanks{
This work was supported in part by the National Natural Science Foundation of China (NSFC) under Grant 52375498,  
and in part by the Fundamental Research Funds for the Central Universities under Grant 21623111. (Corresponding author: Ting Qu)}
\thanks{Wenyou Guo, Hainan Huang and Yafeng Wei are with School of Management, Jinan University, Guangzhou 510632, China (e-mail: guosir1997@163.com; hhn0113@outlook.com; wyf\_gl1022@163.com). }
\thanks{Ting Qu is  with Guangdong International Cooperation Base of Science and Technology for GBA Smart Logistics, Jinan University, Zhuhai 519070, China, also with School of Intelligent Systems Science and Engineering, Jinan University, Zhuhai 519070, China,  and also with Institute of Physical Internet, Jinan University, Zhuhai 519070, China  (e-mail: quting@jnu.edu.cn).}
}


\maketitle

\begin{abstract}
Within the framework of the augmented Lagrangian (AL), we propose a novel distributed optimization method, termed Distributed Augmented Lagrangian Decomposition (DALD), and provide a rigorous convergence proof for its standard version.
To address the high iteration costs in early stages, we propose several accelerated variants  of DALD that enhances efficiency without compromising theoretical guarantees, supported by a comprehensive convergence analysis. To facilitate the description of the distributed optimization process, the concept of hierarchical coordination networks is introduced, integrating hierarchical matrix concepts to aid in this explanation. We further explore and expand the applicability of the DALD method and demonstrate how it unifies existing distributed optimization theories within the AL framework. The effectiveness and applicability of the proposed distributed optimization method and its variants are further validated through numerical experiments.
\end{abstract}

\section{Introduction}
\label{sec:1.introduction}


\IEEEPARstart{I}{n} the current scenario, the decision-making problems faced by various application disciplines are showing a trend of increasing scale and complexity \cite{chatzipanagiotis2015network}. In applications such as transportation, logistics, power systems, smart manufacturing, wireless communications, machine learning, artificial intelligence, and various other fields, these scenarios often come with very large datasets, obtained, stored, and retrieved in a distributed manner. In this context, traditional monolithic optimization methods struggle to effectively address large-scale optimization problems, as data storage and decision-making authority are decentralized, communication bandwidth for information exchange is limited, and privacy concerns for information may arise \cite{giselsson2018large,yangSurveyADMMVariants2022,jeong2023review}. Distributed optimization algorithms that involve multiple computing nodes collaborating can significantly reduce the computational burden on a single machine and are increasingly receiving more attention \cite{fourer2010optimization}.

Research on distributed optimization can be traced back to the pioneering works of Tsitsiklis and Bertsekas \cite{tsitsiklis1984problems,tsitsiklis1986distributed,bertsekas2015parallel}. The core idea of distributed optimization, known as ``decomposition-coordination" \cite{cohen1978optimization}, aims to decompose the all-in-one  (AIO) model into smaller subproblems, each of which can be locally solved by individual computing nodes without needing access to all information of the entire system \cite{berahas2018balancing,nedic2018distributed}. Subsequently, by coordinating these solutions through appropriate mechanisms, the solutions of subproblems gradually converge to the optimal solution of the entire system and satisfy the constraints of the entire system \cite{yfantis2023hierarchical}. This idea effectively exploits the parallel computing capabilities of distributed computer systems, ensures the global objective of the entire system is met, and simultaneously mitigates the risk of data privacy breaches \cite{yuan2023distributed}.

Dual decomposition, as a classical decomposition method, is also known as Lagrangian relaxation or Lagrangian decomposition, which originates from Lagrangian dual theory \cite{lasdon1968duality}. The idea first appeared in \cite{everett1963generalized} and has since been explored in numerous literatures and widely used to solve integer programming problems \cite{geoffrion2009lagrangean,fisher1981lagrangian,guignard2003lagrangean,wolsey2020integer}. The method is computationally simple and popular, but due to the non-differentiability of the dual function induced by the ordinary Lagrangian, which requires the application of advanced non-smooth optimization techniques to ensure numerical stability and efficiency \cite{lee2017complexity}. The Augmented Lagrangian Method  (ALM), also referred to as the ``Method of Multipliers", has become one of the most critical and widely used tools in the field of constraint optimization, primarily due to its superior convergence speed and numerical properties \cite{bertsekas1996constrained,birgin2014practical}.

The specialized techniques for decomposing the AL can be traced back to early works \cite{gabay1976dual,watanabe1978decomposition,glowinski1987augmented,tatjewski1989new}. Recent literatures involve Diagonal Quadratic Approximation (DQA) \cite{mulvey1992diagonal,ruszczynski1995convergence,berger1994extension}, Accelerated Distributed Augmented Lagrangian (ADAL) \cite{chatzipanagiotis2015augmented,chatzipanagiotis2015distributednoise}, and Alternating Direction Method of Multipliers (ADMM) \cite{eckstein1989splitting,eckstein1992douglas}. The core idea of DQA lies in generating separable approximations of the primal step taken in ALM, which iteratively converge to the true primal steps through the addition of correction steps. 
It is worth noting that ADAL, introduced in \cite{chatzipanagiotis2015augmented}, was inspired by DQA and can be viewed as a truncated version of the DQA algorithm. Another approach to decomposition-coordination is ADMM, employing the Gauss-Seidel scheme. Initially, ADMM was primarily used to solve optimization problems with two blocks of variables \cite{boyd2011distributed}, which we refer to as classic ADMM in this paper. Nevertheless, the convergence of ADMM when directly extended to optimization problems involving multiple blocks of variables remains an unresolved issue \cite{chen2016direct}. Furthermore, to address problems in various scenarios, numerous variants of classic ADMM have been proposed by previous researchers \cite{yangSurveyADMMVariants2022}.

However, the aforementioned distributed optimization algorithms typically rely on varying degrees of assumptions, particularly regarding the separability and convexity of the objective function,  and are often limited to a two-level hierarchy in the sequence of subproblem solving. These assumptions are fundamental to their theoretical foundations and practical performance. Moreover, these algorithms typically adhere to a fixed subproblem-solving sequence, either simply serial or parallel, without delving deeply into the internal structure of the problem. Against this backdrop, this paper introduces a novel distributed optimization technique, the Distributed Augmented Lagrangian Decomposition (DALD) method, designed for distributed solutions to more general optimization problems. What sets this method apart is its ability to handle objective functions that may include non-separable terms, while enabling distributed optimization over flexible communication topologies. During the development of this method, we also reveal the essence of existing distributed optimization algorithms based on the AL techniques. 

The main contributions of this paper are as follows:
\begin{itemize}
    \item[a)]  A distributed solver-agnostic optimization technique, DALD, is introduced, along with a rigorous convergence proof for its standard version.

    \item [b)] The concept of hierarchical coordination networks is presented for the first time, systematically characterizing the dynamics of distributed optimization, with hierarchical matrices aiding the formulation.

    \item [c)] Accelerated variants of DALD are developed, achieving significant computational efficiency gains while providing a convergence analysis.

    \item[d)]  Additionally, the paper extends the applicability of DALD to unify existing distributed optimization theories within the AL framework and explains why direct extensions of ADMM may fail to converge~\cite{chen2016direct}, addressing a critical gap in prior theory.
\end{itemize}

The rest of the paper is organized as follows: Section \ref{sec:2.Preliminaries} reviews relevant studies and provides the necessary definitions. Section \ref{sec:3.The Distributed Augmented Lagrangian Decomposition Method} introduces the standard version of the DALD algorithm and its convergence properties. Section \ref{sec:4.DALD in Practice} presents the accelerated version of DALD and outlines a unified framework for distributed optimization. Section \ref{sec:5.Numerical Experiments} details the experimental design and computational results, and Section \ref{sec:6.Conclusion} concludes the paper. To provide readers with a clearer understanding of DALD, the relevant derivations are deferred to Appendix \ref{sec:7.Appendices}.

\section{Preliminaries}
\label{sec:2.Preliminaries}

\subsection{Related Work}
\label{subsection:Related Work}

We first consider the following optimization problem:
\begin{equation}
\begin{aligned}
& \underset{\substack{\mathbf{x} \in X}}{\text{min}}
& & f(\mathbf{x})  \\
&  \text{~s.t.}
& & \sum_{i=1}^n \mathbf{A}_i \mathbf{x}_i = \mathbf{b}
\end{aligned}
\label{eq:1}
\end{equation}
\noindent where the function $f: \mathbb{R}^N \rightarrow \mathbb{R}$ is  continuously differentiable and convex. Furthermore, $\mathbf{x}_i$ are multidimensional subvectors of $\mathbf{x}$, and $\mathbf{x}_i = [x_{i1}, x_{i2}, \ldots, x_{iN_{i}}]^{\top} \in  X_i \subseteq \mathbb{R}^{N_i}$, $i \in \mathbb{N}_n = \{1, 2, \ldots, n\}$. Collectively, $\mathbf{x} = [\mathbf{x}_1^{\top}, \mathbf{x}_2^{\top}, \ldots, \mathbf{x}_n^{\top}]^{\top} \in X \subseteq \mathbb{R}^N$ with $N = \sum_{i=1}^{n} N_i$, and $X$ is a Cartesian product of closed sets: $X = X_1 \times X_2 \times \cdots \times X_n$, where $X_i$ denotes a given nonempty closed subset of $\mathbb{R}^{N_i}$, and each $\mathbf{A}_i$ is a matrix of dimension $m \times N_i$, $\mathbf{b} \in \mathbb{R}^m$  is a vector, $i \in \mathbb{N}_n$.

To facilitate the symbolic representation for subsequent analysis, the constraints in problem (\ref{eq:1}) can be formulated as:
\begin{equation*}
    \varphi(\mathbf{x}) = \sum_{i=1}^n \mathbf{A}_i \mathbf{x}_i - \mathbf{b} =0,
\end{equation*}
where $\varphi: \mathbb{R}^{N} \rightarrow \mathbb{R}^m$.
Subsequently, by introducing Lagrange multipliers \( \mu \in \mathbb{R}^m \) and penalty parameters \( \rho \in \mathbb{R}^m \), we relax problem \eqref{eq:1} and define the Lagrangian as:  
\begin{equation}
    L(\mathbf{x}, \mu) = f(\mathbf{x}) + \mu^{\top} \varphi(\mathbf{x}), \label{eq:2}
\end{equation}
and the AL is defined as:
\begin{equation}
    \Lambda_{\rho}(\mathbf{x}, \mu) = f(\mathbf{x}) + \mu^{\top} \varphi(\mathbf{x}) + \|\rho \circ \varphi(\mathbf{x})\|^2. \label{eq:3}
\end{equation}
The symbol $\circ$ represents the Hadamard product: $\mathbf{c} \circ \mathbf{d}  =  [c_1, c_2, \ldots, c_n]^{\top} \circ [d_1, d_2, \ldots, d_n]^{\top} = [c_1 d_1, c_2 d_2, \ldots, c_n d_n]^{\top}$, and $\|\cdot\|$ denotes the $\ell_2$-norm.

To tackle the problem \eqref{eq:1}, we can leverage the ALM.
\begin{algorithm}[H]
\caption{Augmented Lagrangian Method (ALM)}\label{alg:alm}
\begin{algorithmic}
\STATE 
 \textbf{Initialize}: Set $k = 1$, and give the initial Lagrange multipliers $\mu^1$ and the penalty $\rho$.
\STATE 
\textbf{Step 1}: For fixed parameters $\mu^k$, calculate $\mathbf{x}^k$ as a solution of the problem:
\begin{equation}
    \mathbf{x}^k = \arg \min_{\mathbf{x} \in X} \Lambda_{\rho} (\mathbf{x}, \mu^k). \label{eq:4}
\end{equation}
\vspace{-1em}
\STATE 
 \textbf{Step 2}: If $\varphi(\mathbf{x}^k) = 0$, then stop (optimal solution found), let $\mathbf{x}^* = \mathbf{x}^k$. Otherwise, update:
\begin{equation}
    \mu^{k+1} = \mu^k + 2\rho \circ \rho \circ\varphi(\mathbf{x}^k). \label{eq:5}
\end{equation}
\hspace{2em} 
Set $k = k + 1$, and repeat from Step 1.
\end{algorithmic}
\end{algorithm}

The ALM entails a two-layer process: firstly, an inner layer  where the relaxation problem is minimized--given the current Lagrange multipliers $\mu^k$, the optimal solution $\mathbf{x}^k$ is derived by solving $\Lambda_\rho (\mathbf{x}, \mu^k)$. Secondly, an outer layer is performed, focusing on updating the Lagrange multipliers from $\mu^k$ to $\mu^{k+1}$ to refine the adherence to the problem's constraints. Under the stipulations outlined in assumptions \hyperref[assumption:A1]{(A1)}–\hyperref[assumption:A3]{(A3)}, 
which will be detailed in Section~\ref{sec:3.2 Convergence Analysis of Standard DALD}, 
employing ALM for tackling the problem (\ref{eq:1}) yields the optimal solution~\cite{bertsekasMethodMultipliersConvex1975,bertsekas_nonlinear_2016,bertsekas1996constrained,xuIterationComplexityInexact2021}.

\subsection{Notation}
\label{subsection:Notation}

\begin{figure}[t!bp]
\centering
\includegraphics[scale=0.32]{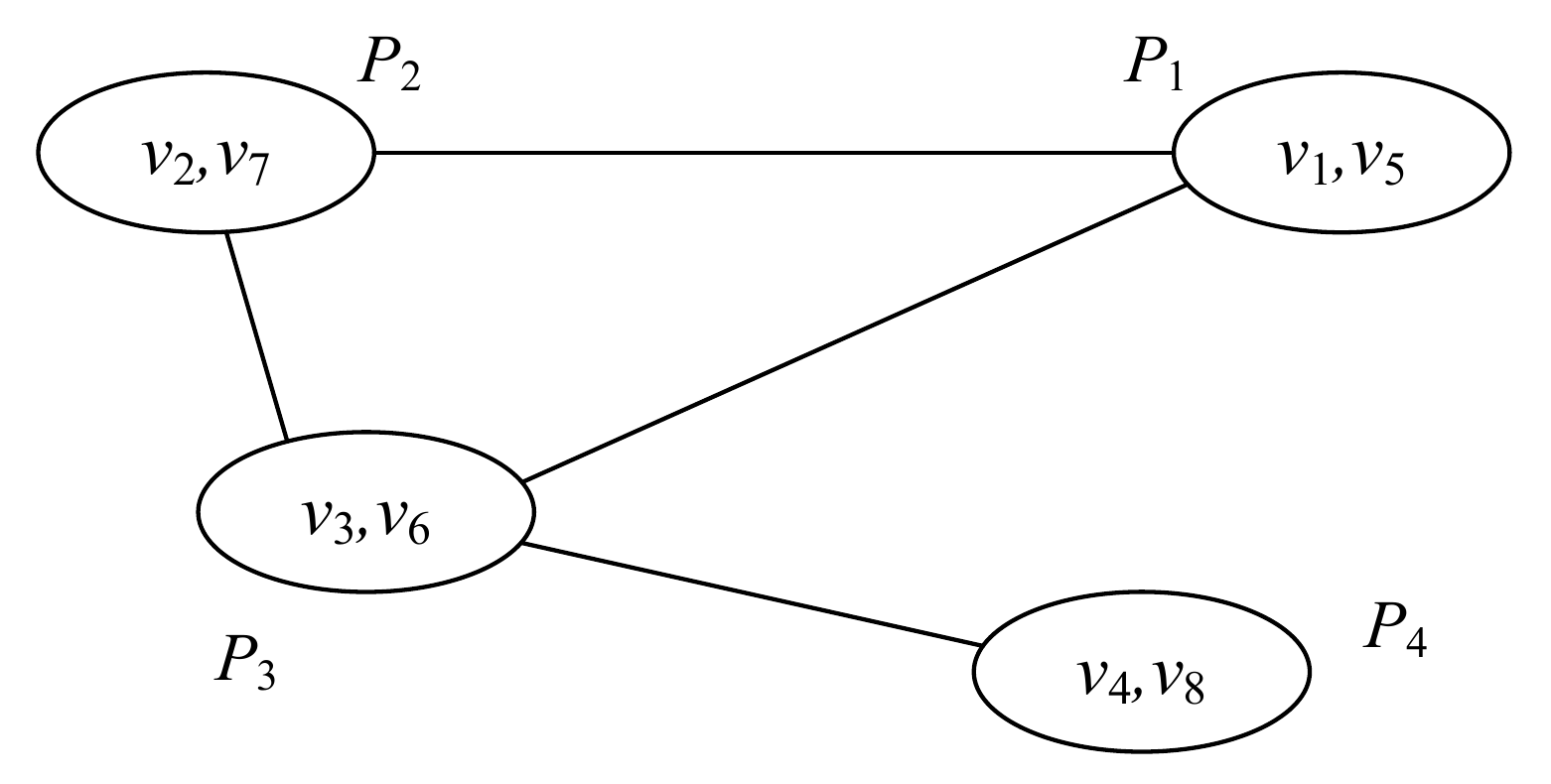}
\caption{Distribution of decision-making variables and coupling relationships}
\label{fig:1}
\end{figure}

Let us consider a large-scale optimization problem \( P \) formulated as (\ref{eq:1}), which can be decomposed into \( n \) subproblems \( P_i \), \( i \in \mathbb{N}_n \). 
In the following, we introduce some of the concepts and definitions adopted in this paper. For clarity, we use a simple example with \( n = 4 \) to illustrate the distribution of decision variables and the coupling relationships among these  subproblems, as shown in Fig.~\ref{fig:1}. 
The optimization problem is defined as follows:
\begin{equation*}
\begin{aligned}
& \min
& &  \underbrace{ v_3 + v_4 + v_5^2 + v_7^2 + v_6^2 + v_8^2}_{\sum_{i=1}^{n} f_i(\mathbf{x}_i)}
\;+\;
\underbrace{2v_5 v_7 + 2v_6 v_8}_{f_0(\mathbf{x}_1,\mathbf{x}_2,\ldots,\mathbf{x}_n)}\\
& ~\text{s.t.} 
& &  \varphi(\mathbf{x}) = \left\{
\begin{aligned}
& v_1 + v_5 + v_6 - 7=0,  \\
& v_2 + v_5 + v_7 - 10=0, \\
& v_5+v_6 + 2v_7 - 12=0,
\end{aligned}
\right. \\
& & & \mathbf{x} \in X:~  1 \leq v_i \leq 7,  i \in \mathbb{N}_8.
\end{aligned}
\end{equation*}

\renewcommand{\labelitemi}{$\bullet$}
\begin{itemize}
    \item If the variable $\mathbf{x}_i$ exists in subproblem $P_i$, then $\mathbf{x}_i$ is designated as the \textit{local variable} of subproblem $P_i$, and the local variables between subproblems do not intersect.
\end{itemize}

Unlike the notation used in Section \ref{subsection:Related Work}, we denote the vector $\mathbf{x}$ as $\mathbf{x} = (\mathbf{x}_1, \mathbf{x}_2, \ldots, \mathbf{x}_n)$. Similar adjustments were made to other related vector symbols, and we slightly abused some of them. 
For the above problem, we set $\mathbf{x}_1 = (v_1, v_5)$, $\mathbf{x}_2 = (v_2, v_7)$, $\mathbf{x}_3 = (v_3, v_6)$, $\mathbf{x}_4 = (v_4, v_8)$.

\textbf{Definition 1.} The optimization problem $P$ formulated as (\ref{eq:1}) can be decomposed into $n$ subproblems. Subproblems $P_a$ and $P_b$ are \textit{coupled} if:
\begin{enumerate}
    \item Certain non-separable terms in the objective function of problem $P$ depend on elements from  both $\mathbf{x}_a$ and $\mathbf{x}_b$,
    \item or some elements of  $\mathbf{x}_a$ and $\mathbf{x}_b$ appear in the same constraint, $a, b \in \mathbb{N}_n$, $a \neq b$.
\end{enumerate}

\renewcommand{\labelitemi}{$\bullet$}
\begin{itemize}
    
    \item If subproblem $P_i$ is coupled with subproblem $P_j$, then the set $\mathcal{R}_i$ represents the \textit{coupling relationships} of all subproblems $P_j$ that are associated with subproblem $P_i$, i.e., $j(P_j) \in \mathcal{R}_i$, and the following relationships exist: $\mathcal{R}_1 = \{2, 3\}, \; \mathcal{R}_2 = \{1, 3\}, \; \mathcal{R}_3 = \{1, 2, 4\}, \; \mathcal{R}_4 = \{3\}.$

    \item $\mathbf{x}_{-i}$ represents the \textit{coupling variable} of the subproblem $P_i$, where $\mathbf{x}_{-i} = (\mathbf{x}_j ^ {S_i}\mid j \in \mathcal{R}_i)$, and $\mathbf{x}_j ^ {S_i}$  denotes the subset of elements from $\mathbf{x}_j$ for $\mathbf{x}_{-i}$, $i \in \mathbb{N}_n$. 
    For $i \neq j$, there may be common elements between $\mathbf{x}_{-i}$ and $\mathbf{x}_{-j}$.
     We have $\mathbf{x}_{(-1)} = (v_6, v_7)$, $\mathbf{x}_{(-2)} = (v_5, v_6)$, $\mathbf{x}_{(-3)} = (v_5, v_7, v_8)$, $\mathbf{x}_{(-4)} = (v_6)$.

    \item 
    In the objective function, we define $f_0(\mathbf{x})$ as the cumulative sum of all non-separable sub-terms across different subproblems. We further define $f_0(\mathbf{x}_i, \mathbf{x}_{-i})$ as the sum of the non-separable terms in $f_0(\mathbf{x})$ that are associated with subproblem $i$ and its coupling subproblems, which is referred to as the \textit{coupling objective function}.  Let $f_i(\mathbf{x}_i)$ be the \textit{local objective function} for subproblem $P_i$, and $f_i(\mathbf{x}_i, \mathbf{x}_{-i}) = f_i(\mathbf{x}_i) + f_0(\mathbf{x}_i, \mathbf{x}_{-i})$, $i \in \mathbb{N}_n$. Note that $f_0(\mathbf{x})$ is not equivalent to $\sum f_0(\mathbf{x}_i, \mathbf{x}_{-i})$. To illustrate, we have $f_1(\mathbf{x}_1) = v_5^2, f_0(\mathbf{x}_1, \mathbf{x}_{-1}) = 2v_5v_7,  f_1(\mathbf{x}_1, \mathbf{x}_{-1}) =  v_5^2 + 2v_5v_7,  f_0(\mathbf{x}) = 2v_5 v_7 + 2v_6 v_8.$

    \item Constraint $\phi_i(\mathbf{x}_i, \mathbf{x}_{-i}) = 0$, which contains the local variable $\mathbf{x}_i$ and $\mathbf{x}_j$, is called the \textit{coupling constraint} of subproblem $P_i$ and subproblem $P_j$, $i \in \mathbb{N}_n$, $j \in \mathcal{R}_i$.  As clearly indicated in 2) of Definition 1, when subproblems are coupled, their coupling constraints may share common elements.   For example, $ \phi_2(\mathbf{x}_2, \mathbf{x}_{-2}) $ consists of the elements in the set $ \{v_2 + v_5 + v_7 - 10 =0, \, v_5 + v_6 + 2v_7 - 12=0\} $ and $ \phi_3(\mathbf{x}_3, \mathbf{x}_{-3}) $ consists of the elements in the set $ \{v_1 + v_5 + v_6 - 7=0, \,  v_5 + v_6 + 2v_7 - 12=0\}. $

    \item The local variable $\mathbf{x}_i$ is drawn from the \textit{local constraint} set $X_i$, i.e., $\mathbf{x}_i \in X_i$, $i \in \mathbb{N}_n$. For instance, the local constraint set $ \mathbf{x}_3 $ is defined by the numerical interval $ 1 \leq \mathbf{x}_3 \leq 7 $.
\end{itemize}

In the expressions $f_0(\mathbf{x}_i, \mathbf{x}_{-i})$ and $\phi_i(\mathbf{x}_i, \mathbf{x}_{-i})$, the coefficients and operations reflect the shared knowledge between subproblem (subsystem or agent) $P_i$  and other subproblems $P_j$, $i \in \mathbb{N}_n$, $j \in \mathcal{R}_i$. 
All relevant subproblems in $\{i\} \cup \mathcal{R}_i$ are generally considered to have access to this shared information. In contrast, $f_i(\mathbf{x}_i)$ and $X_i$ are considered private knowledge of $P_i$ and should not be disclosed to the outside, $i \in \mathbb{N}_n$.

\section{The Distributed Augmented Lagrangian Decomposition Method}
\label{sec:3.The Distributed Augmented Lagrangian Decomposition Method}

\subsection{Standard Version}
\label{sec:3.1Standard Version}
Building upon the notation established in Section \ref{subsection:Notation}, we reformulate problem (\ref{eq:1}) as follows:
\begin{equation}
\begin{aligned}
& \underset{\substack{\mathbf{x} \in X}}{\min}
& & \sum_{i=1}^{n} f_i (\mathbf{x}_i) + f_0 (\mathbf{x}_1, \mathbf{x}_2, \ldots, \mathbf{x}_n) \\
& \text{~s.t.}
& & \phi_i (\mathbf{x}_i, \mathbf{x}_{-i}) = 0,  i \in \mathbb{N}_n,
\end{aligned}
\label{eq:6}
\end{equation}
where $\phi_i : \mathbb{R}^{N} \rightarrow \mathbb{R}^{m_i}$, $i \in \mathbb{N}_n$. The vector $\varphi(\mathbf{x})$, introduced in~\ref{subsection:Related Work},
is constructed by concatenating the elements of all vectors $\{\phi_i (\mathbf{x}_i, \mathbf{x}_{-i}) \mid i \in \mathbb{N}_n\}$, with shared elements included only once. Subsequently, problem (\ref{eq:6}) is decomposed into individual subproblems $P_i$:
\begin{equation}
\begin{aligned}
& \underset{\substack{\mathbf{x}_i \in X_i}}{\min}
& & f_i (\mathbf{x}_i, \mathbf{x}_{-i}) \\
& \text{~~s.t.}
& & \phi_i (\mathbf{x}_i, \mathbf{x}_{-i}) = 0.
\end{aligned}
\label{eq:7}
\end{equation}

We introduce local Lagrange multipliers $\mu_i \in \mathbb{R}^{m_i}$ and local penalty parameters $\rho_i \in \mathbb{R}^{m_i}$ to relax problem (\ref{eq:7}), $i \in \mathbb{N}_n$. The local Lagrangian takes the form:
\begin{equation}
\mathcal{L}_i (\mathbf{x}_i, \mathbf{x}_{-i}, \mu_i) = f_i (\mathbf{x}_i, \mathbf{x}_{-i}) + \mu_i^{\top} \phi_i (\mathbf{x}_i, \mathbf{x}_{-i}).
\label{eq:8}
\end{equation}

Furthermore, the local AL is defined as:
\begin{equation}
\begin{aligned}
\mathcal{A}_{\rho_i}^i (\mathbf{x}_i, \mathbf{x}_{-i}, \mu_i) = & {} f_i (\mathbf{x}_i, \mathbf{x}_{-i}) + \mu_i^{\top} \phi_i (\mathbf{x}_i, \mathbf{x}_{-i}) \\
 & {} +  \| \rho_i \circ \phi_i (\mathbf{x}_i, \mathbf{x}_{-i}) \|^2.
\end{aligned}
\label{eq:9}
\end{equation}
Consequently,  the vectors $\mu$ and $\rho$, mentioned in Section \ref{subsection:Related Work}, are constructed by concatenating the unique elements of all vectors  $\{\mu_i \mid i \in \mathbb{N}_n\}$ and $\{\rho_i \mid i \in \mathbb{N}_n\}$, respectively.

In our application of ALM to complex large-scale problems, we acknowledge the interdependencies among subproblems and decompose the primal problem (\ref{eq:6}) into subproblems, optimizing them alternatively with the Block Coordinate Descent (BCD) algorithm idea. Stemming from this idea, the DALD algorithm emerges, comprising two stages: initially, within the inner layer loop, the relaxation problem is decomposed and minimized with $\mathbf{x}_1, \mathbf{x}_2, \ldots, \mathbf{x}_n$ as variable components, yielding the optimal solution $\mathbf{x}^{k,*}$ for a given $\mu^k$. Then, the outer layer loop updates the Lagrange multipliers to $\mu^{k+1}$ based on $\mathbf{x}^{k,*}$.

We define the notation \( \mathbf{x}_i^{k,v} \) to represent the local variable computed by subproblem \( i \) during the \( v \)-th inner iteration of the \( k \)-th outer loop, i.e., the \( (k,v) \)-th iteration of the DALD.
For notational simplicity,
 we introduce a new set of auxiliary variables defined as follows for $i \in \mathbb{N}_n$:
\begin{equation}
\begin{aligned}
\mathbf{w}_i^{k,v} &= (\mathbf{x}_1^{k,v}, \ldots, \mathbf{x}_{i-1}^{k,v}, \mathbf{x}_i^{k,v}, \mathbf{x}_{i+1}^{k,v-1}, \ldots, \mathbf{x}_n^{k,v-1}). 
\end{aligned}
\label{eq:10}
\end{equation}

Moreover,
\begin{equation}
\begin{aligned}
\mathbf{w}_{-i}^{k,v} = \{ \mathbf{x}_j^{k, v} \mid j \in \mathcal{R}_i, \, j < i \} \cup  \{ \mathbf{x}_j^{k, v-1} \mid j \in \mathcal{R}_i, \, j > i \}.
\end{aligned}
\label{eq:11}
\end{equation}

Next, we present the standard DALD algorithm, and its working principle is as follows:

\begin{algorithm}[H]
\caption{Standard Distributed Augmented Lagrangian Decomposition (Standard DALD)}\label{alg:dald}
\begin{algorithmic}
\STATE 
\textbf{Initialize}: Set $k = 1$, $v = 0$, and choose $\mathbf{x}^{1,0}$, $\mu^{1}$, $\rho$.
\STATE 
 \textbf{Step 1.1}: Let $v = v + 1$. Given a sequence $\mathbb{S}$, perform sequence $\mathbb{S}$ to solve subproblems $P_s$, $s \in \mathbb{S}$:
\begin{equation}
    \mathbf{x}_s^{k,v} = \arg \min_{\mathbf{x}_s \in X_s} \mathcal{A}_{\rho_s}^s (\mathbf{x}_s, \mathbf{w}_{-s}^{k,v}, \mu_s^k). \label{eq:12}
\end{equation}
\vspace{-1em}
\STATE 
\textbf{Step 1.2}: Until the last subproblem is solved, resulting in $\mathbf{x}^{k,v}$, 
$\mathcal{C}_i^{k,v}= \phi_i (\mathbf{x}_i^{k,v}, \mathbf{x}_{-i}^{k,v})$, and $\mathcal{D}_i^{k,v}= \mathbf{x}_i^{k,v} - \mathbf{x}_i^{k,v-1}, i \in \mathbb{N}_n$.
If  $\mathcal{D}^{k,v} = 0$ is satisfied, proceed to Step 2; otherwise, return to Step 1.1.
\STATE 
\textbf{Step 2}: If $\mathcal{C}^{k,v} = 0$ is satisfied, then stop (optimal solution found), let $\mathbf{x}^* = \mathbf{x}^{k,v}$. Otherwise, update:
\begin{equation}
    \mu_i^{k+1} = \mu_i^k + 2 \rho_i \circ \rho_i \circ \mathcal{C}_i^{k,v},  i \in \mathbb{N}_n. \label{eq:13}
\end{equation}
\hspace{2em} 
Set $\mathbf{x}^{k+1,0} = \mathbf{x}^{k,v}$, $k = k + 1$, $v = 0$, and repeat from Step 1.1.
\end{algorithmic}
\end{algorithm}

Here, $\mathcal{C}^{k,v}$ and $\mathcal{D}^{k,v}$ denote the primal residual and dual residual at the $(k,v$)-th iteration of the DALD algorithm, respectively. Detailed  definitions and calculation procedures  are provided in Appendix \ref{subsection:appendix A}.

For our subsequent discussions, we assume that $\{\mathbf{x}^{k,v}\}$ is the sequence generated during the $k$-th outer loop iteration of the DALD. The limit point or endpoint of this sequence is denoted by $\bar{\mathbf{x}}^{k,v}$, and $\bar{\mathbf{x}}_i^{k,v-1}$ represents the element immediately preceding $\bar{\mathbf{x}}_i^{k,v}$. The sequence $\{\mathbf{x}^k\}$, derived from $\bar{\mathbf{x}}^{k,v}$, has a limit point $\mathbf{x}^*$. Similarly, $\bar{\mu}^k$, $\mu^*$, and $\{\mu^k\}$ are obtained.

\subsection{Convergence Analysis of Standard DALD}
\label{sec:3.2 Convergence Analysis of Standard DALD}

In this section, we provide a theoretical convergence analysis of the standard version of the DALD algorithm. Based on the discussion in Section~\ref{sec:3.1Standard Version} regarding the design motivation of the DALD algorithm, the proof of its convergence can be simplified to verifying the convergence of its inner layer. Specifically, for a given set of Lagrange multipliers \(\mu^k\), it is necessary to demonstrate that the inner layer iterations converge to \(\mathbf{x}^{k,*}\) within the \(k\)-th outer loop.  Once this is verified, the remaining task in the outer layer loop is simply to update the Lagrange multipliers.

Before we proceed, we first need to make the following assumptions regarding problem (\ref{eq:1}):

\begin{enumerate}
\item[(A1)] The function $f(\mathbf{x})$ is continuously differentiable and convex. \label{assumption:A1}
\end{enumerate}

\begin{enumerate}
\item[(A2)] Each $X_i \subseteq \mathbb{R}^{N_i}$ is a nonempty, closed, and convex subset, $i \in \mathbb{N}_n$. \label{assumption:A2}
\end{enumerate}

\begin{enumerate}
\item[(A3)] The Lagrangian (\ref{eq:2}) has a saddle point $(\mathbf{x}^*, \mu^*) \in \mathbb{R}^N \times \mathbb{R}^m$:
\label{assumption:A3}
\end{enumerate}
\begin{equation*}
    L(\mathbf{x}^*, \mu) \leq L(\mathbf{x}^*, \mu^*) \leq L(\mathbf{x}, \mu^*), \forall \mathbf{x} \in X, \forall \mu \in \mathbb{R}^m.
\end{equation*}

\begin{enumerate}
\item[(A4)] All  subproblems are solvable at each iteration. 
\label{assumption:A4}
\end{enumerate}

Assumptions \hyperref[assumption:A1]{(A1)}–\hyperref[assumption:A3]{(A3)} imply that (\ref{eq:4}) and (\ref{eq:12}) are convex, and the strong duality holds for the problem (\ref{eq:1}), i.e., the optimal values of the primal and dual problems are equal.
Under these assumptions, since all constraints $\varphi(\mathbf{x})$ are affine functions, and when a feasible solution exists, the Slater condition is automatically satisfied, thereby guaranteeing that \hyperref[assumption:A3]{(A3)} holds.
Furthermore, these assumptions guarantee the convergence of the ALM. Assumption \hyperref[assumption:A4]{(A4)} ensures the smooth implementation of the inner loops of DALD.

\begin{figure*}[t!b]
\centering
\includegraphics[scale=0.46]{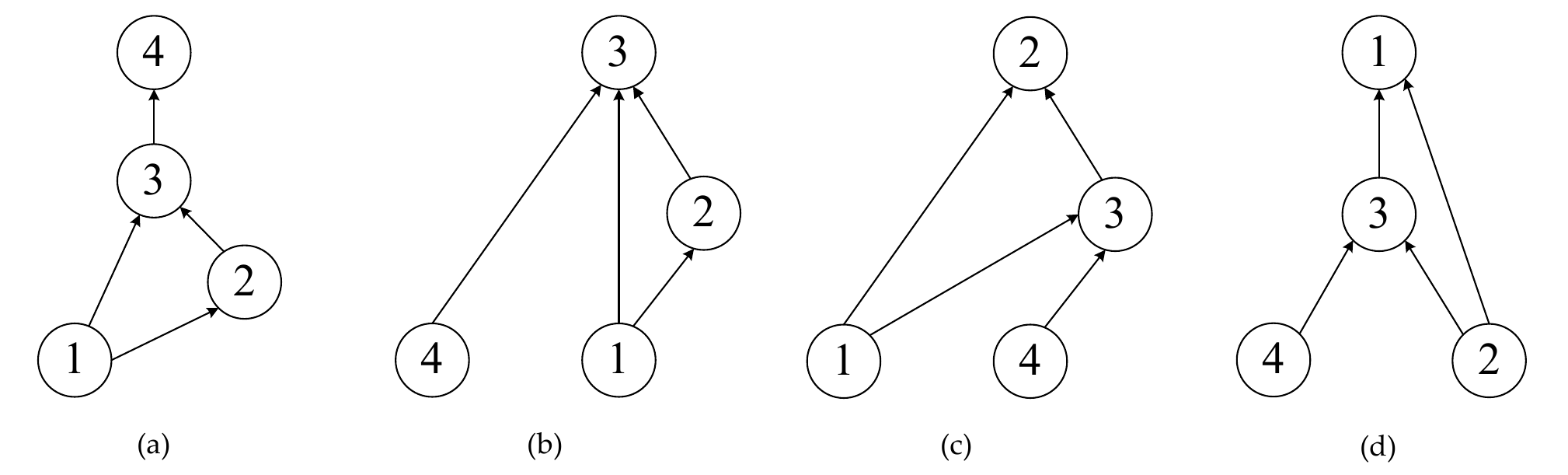}
\caption{ Hierarchical networks for describing the subproblem solving sequences }
\label{fig:2}
\end{figure*}

\begin{figure*}[t!b]
\centering
\includegraphics[scale=0.50]{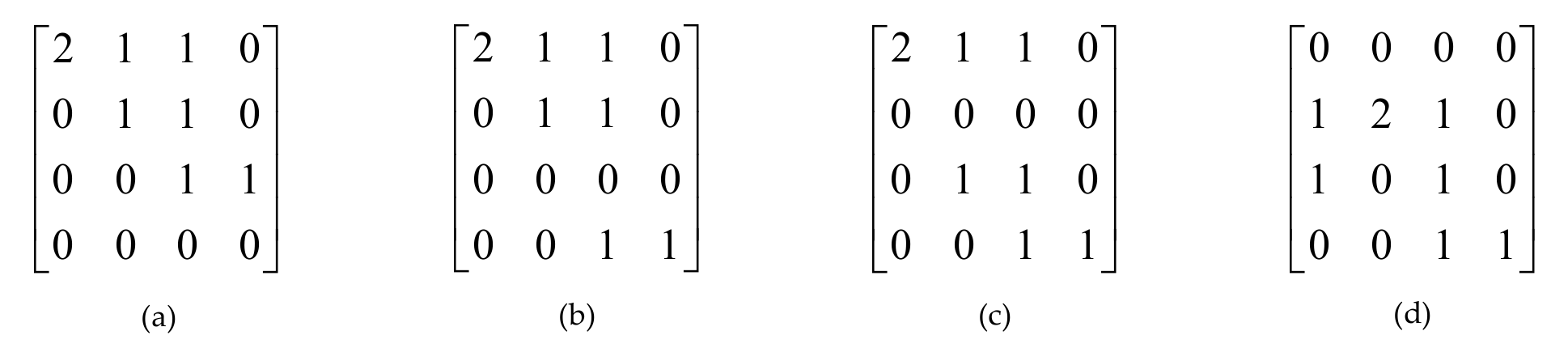}
\caption{ Hierarchical matrices for describing the subproblem solving sequences }
\label{fig:3}
\end{figure*}


The subsequent theorem establishes the convergence of the algorithm's inner layer loops, thereby proving the convergence of the entire algorithm.

\noindent \textbf{Theorem 1} Assume \hyperref[assumption:A1]{(A1)}--\hyperref[assumption:A2]{(A2)} and  \hyperref[assumption:A4]{(A4)} . Then, every limit point of $\{\mathbf{x}^{k,v}\}$ minimizes $\Lambda_\rho (\mathbf{x}, \mu^k)$ over $X$.

\textbf{\textit{Proof:}}
When $\mu^k$ and $\rho$ are fixed, we denote
\begin{equation}
    \mathcal{F}(\mathbf{x}) = \Lambda_\rho (\mathbf{x}, \mu^k). \label{eq:14}
\end{equation}
By assumption \hyperref[assumption:A1]{(A1)}, $\mathcal{F}(\mathbf{x})$ is convex.
 Under assumption \hyperref[assumption:A4]{(A4)} , given the current iterate $\mathbf{x}^{k,v} = (\mathbf{x}_1^{k,v}, \mathbf{x}_2^{k,v}, \ldots, \mathbf{x}_n^{k,v})$, the next iterate $\mathbf{x}^{k,v+1} = (\mathbf{x}_1^{k,v+1}, \mathbf{x}_2^{k,v+1}, \ldots, \mathbf{x}_n^{k,v+1})$ is computed via (12). By optimizing $\mathbf{x}_i$ sequentially, we derive
\begin{equation}
    \mathcal{F}(\mathbf{x}^{k,v}) \geq \mathcal{F}(\mathbf{w}_1^{k,v+1}) 
    \geq \cdots \geq \mathcal{F}(\mathbf{w}_{n-1}^{k,v+1}) \geq \mathcal{F}(\mathbf{x}^{k,v+1}). \label{eq:15}
\end{equation}
Let $\bar{\mathbf{x}}^{k,v} = (\bar{\mathbf{x}}_1^{k,v}, \bar{\mathbf{x}}_2^{k,v}, \ldots, \bar{\mathbf{x}}_n^{k,v})$ be a limit point of the sequence $\{\mathbf{x}^{k,v}\}$. Equation (\ref{eq:15}) implies that the sequence $\{\mathcal{F}(\mathbf{x}^{k,v})\}$ converges to $\mathcal{F}(\bar{\mathbf{x}}^{k,v})$. We will now demonstrate that $\bar{\mathbf{x}}^{k,v}$ satisfies the optimality condition
\begin{equation*}
    \nabla \mathcal{F}(\bar{\mathbf{x}}^{k,v})^{\top} (\mathbf{x} - \bar{\mathbf{x}}^{k,v}) \geq 0, \forall \mathbf{x} \in X.
\end{equation*}
Consider the sequence $\{\mathbf{x}^{k,v}\}$ converges to $\bar{\mathbf{x}}^{k,v}$. From (\ref{eq:12}) of the DALD and (\ref{eq:15}), we deduce
\[
    \mathcal{F}(\mathbf{x}^{k,v+1}) \leq \mathcal{F}(\mathbf{w}_1^{k,v}) \leq \mathcal{F}(\mathbf{x}_1, \mathbf{x}_2^{k,v}, \ldots, \mathbf{x}_n^{k,v}), \forall \mathbf{x}_1 \in X_1.
\]
Taking the limit as $v \to \infty$, we  obtain
\[
    \mathcal{F}(\bar{\mathbf{x}}^{k,v}) \leq \mathcal{F}(\mathbf{x}_1, \bar{\mathbf{x}}_2^{k,v}, \ldots, \bar{\mathbf{x}}_n^{k,v}), \forall \mathbf{x}_1 \in X_1.
\]
Since $\mathcal{F}$ is convex, its restriction to any coordinate (with the other variables fixed) is also convex.
By the first-order optimality condition (Prop. 3.1.1 in \cite{bertsekas_nonlinear_2016}), it follows that
\begin{equation*}
    \nabla_1 \mathcal{F}(\bar{\mathbf{x}}^{k,v})^{\top} (\mathbf{x}_1 - \bar{\mathbf{x}}_1^{k,v}) \geq 0, \forall \mathbf{x}_1 \in X_1,
\end{equation*}
where $\nabla_i$ denotes the gradient of $\mathcal{F}$  with respect to $\mathbf{x}_i$.

Similarly, for each component $\mathbf{x}_i$
\begin{equation}
    \nabla_i \mathcal{F}(\bar{\mathbf{x}}^{k,v})^{\top} (\mathbf{x}_i - \bar{\mathbf{x}}_i^{k,v}) \geq 0, \forall \mathbf{x}_i \in X_i, i \in \mathbb{N}_n. \label{eq:16}
\end{equation}
Summing inequalities (\ref{eq:16}) and utilizing the Cartesian product structure of the set $X$, we deduce
\begin{equation*}
    \nabla \mathcal{F}(\bar{\mathbf{x}}^{k,v})^{\top} (\mathbf{x} - \bar{\mathbf{x}}^{k,v}) \geq 0, \forall \mathbf{x} \in X.
\end{equation*}
This completes the proof.\hfill $\square$


\textit{Remark 1:} For problem (\ref{eq:1}), suppose that $\varphi(\mathbf{x})$ represents general constraints (not limited to affine functions). Then either the AL in (\ref{eq:3}) remains convex, in which case Theorem 1 still holds, or every limit point of $\{\mathbf{x}^{k,v}\}$ is a critical point of $\Lambda_{\rho}(\mathbf{x}, \mu^k)$ over $X$.

\section{DALD in Practice}
\label{sec:4.DALD in Practice}
\subsection{Topological Perspective}

Based on inherent assumptions, we typically default to a sequential solving sequence $\mathbb{S}$ in the DALD algorithm, where subproblems are tackled one after another. For instance, as depicted in Fig.~\ref{fig:2}(a), subproblems ``1$\rightarrow$2$\rightarrow$3$\rightarrow$4" are coordinated and solved sequentially. Considering the coupling relationships among subproblems, we propose a more robust and flexible subproblem coordination strategy. To describe the subproblem solving sequence $\mathbb{S}$, we introduce the concepts of Hierarchical Network and Hierarchical Matrix.

\textbf{Definition 2.} A \textit{hierarchical network} is constructed by starting from a designated root node and recursively branching out to lower-level nodes in a top-down fashion. Lower-level nodes are methodically endowed with both a level identifier and a unique node identifier. This process begins with nodes that have already been assigned an identifier and iterates until all nodes are integrated into the hierarchy.

When a hierarchical network has $h$ levels, we refer to the sequence of subproblem solving as an $h$-level hierarchical coordination network, and we call this process $h$-level coordination for short. In a hierarchical network, the root node represents the final subproblem to be addressed, while nodes at lower levels correspond to subproblems that must be resolved before those at higher levels. Directed edges in the network represent information flow from lower-level subproblems to higher-level ones. A parent subproblem begins computation only after all dependent lower-level subproblems have completed and transmitted their information. Subproblems at the same level lack directed arc connections and can be solved in parallel. However, their coupling must be addressed when designing the solving sequence in hierarchical networks due to explicit information pathways. Fig.~\ref{fig:2} illustrates some solving sequences of subproblems through hierarchical networks.

To shorten the total computation time and fully account for the synchronous nature of distributed computing, we propose designing a hierarchical network based on the ``Earliest Start (ES)'' rule. As illustrated in Fig.~\ref{fig:2}(b), since both the computation of subproblems and the transmission of information incur time costs, our goal is to initiate the computation of subproblem $P_4$ as early as possible.

\textbf{Definition 3.}
A \textit{hierarchical matrix} \( \mathcal{H} = [a_{ij}] \) is an \( n \times n \) square matrix that encodes the hierarchical relationships among nodes in a network. Specifically, if node \( i \) is a direct descendant of node \( j \) (i.e., there is a directed edge from \( i \) to \( j \) with no intermediate nodes), then \( a_{ij} = 1 \); otherwise, \( a_{ij} = 0 \), which includes the absence of a connection or the presence of a reverse connection from \( j \) to \( i \). The diagonal entry \( a_{ii} \) represents the out-degree of node \( i \).

As shown in Fig.~\ref{fig:3}, the hierarchical matrices provide a detailed representation of the subproblems solving sequences outlined in Fig.~\ref{fig:2}, thereby enhancing our understanding of how hierarchical networks can be used to represent the subproblem coordination process.

In the aforementioned examples (including DLAD-CC), our default subproblem coordination process adopts \textit{full-cycle coordination}, involving the handling of each subproblem once per iteration within the inner loop. However, to meet the demands for more refined control over the subproblems solving process, our algorithm also supports specialized coordination types, such as \textit{partial-cycle coordination} and \textit{selective-repetitive coordination}. The former allows for the resolution of only a portion or even a single subproblem per iteration in the inner loop, whereas the latter permits certain subproblems to be revisited and solved multiple times, while others may be solved only once or skipped.

The selection of subproblems to be solved can be based on either a random or greedy strategy. However, it is important to emphasize that during the DALD iteration process, each subproblem should be considered with equal probability for potential optimization, though the actual selection may vary. This ensures that no subproblem is over-sampled or under-sampled,  thereby preventing convergence bias and ultimately avoiding the risk of the algorithm converging to an incorrect solution.

Furthermore, the DALD algorithm is capable of configuring varying subproblem solving sequences across different inner loops. Nevertheless, within the scope of this study, we assume that under full-cycle coordination, the subproblem solving sequence remains consistent across each iteration of the inner loop. This assumption helps us to simplify the analysis and algorithm implementation.

\subsection{Parameters Determination and Accelerated Version}
\label{subsection:4.2 Parameters Determination and Accelerated Version}

When applying DALD in practice, it is important to leverage any prior knowledge so that the initial multiplier $\mu^1$ is chosen as close as possible to the optimal multiplier $\mu^*$. The penalty parameter $\rho$ should also be set with care, since an excessively large value may lead to ill-conditioning. In addition, appropriate stopping conditions need to be specified for the inner and outer loops, typically with different termination tolerances: $\epsilon_{\text{dual}}$ for the inner loop and $\epsilon_{\text{pri}}$ for the outer loop. These tolerances define the predetermined acceptable levels of computational accuracy.

In the DALD implementation, the inner loop terminates when the following condition is satisfied:
\begin{equation*}
\label{condition:B1}
    \text{(B1)} ~\|\mathcal{D}^{k,v} \|_\infty \leq \epsilon_{\text{dual}} \approx 0.
\end{equation*}
This ensures that the outer loop receives an accurate solution. Here, $\|\mathcal{D}^{k,v} \|_\infty=\max_{i \in \mathbb{N}_n} \| \mathcal{D}_i^{k,v} \|_\infty$, and $\| \mathcal{D}_i^{k,v} \|_\infty = \max_{j \in \mathbb{N}_{N_i}} |\mathcal{D}_{ij}^{k,v}|$. However, empirical evidence demonstrates that setting excessively high precision in the early stages is a gratuitous waste of computational resources.

The augmented Lagrangian is conventionally minimized exactly; however, the algorithm can still achieve convergence even when the minimization process is prematurely terminated~\cite{bertsekas1975combined}. This observation prompts us to reconsider the necessity of exactly solving $\mathbf{x}^{k,*}$ at the initial phase of DALD. Therefore, we introduce a modified stopping criterion in Step 1.2 to enhance efficiency:

\begin{enumerate}
\label{condition:B2}
\item[(B2)] If the condition $\|\mathcal{D}^{k,v} \|_\infty \leq \epsilon_{\text{dual}}^k$ is satisfied, where $\{\epsilon_{\text{dual}}^k\}$ satisfies $0 \leq \epsilon_{\text{dual}}^k, \forall k$, and $\epsilon_{\text{dual}}^k \rightarrow \epsilon_{\text{dual}}$, proceed to Step 2.
\end{enumerate}

The following lemma is paramount to the convergence analysis, with a particular focus on the criterion for terminating the inner loop, designated as condition \hyperref[condition:B2]{(B2)}.

\noindent \textbf{Lemma 1} Suppose $\mathbf{x}^*$ is the optimal point of problem (\ref{eq:1}). Then, there exists a constant vector $\gamma = (\gamma_1, \gamma_2, \ldots, \gamma_n) \in \mathbb{R}^N$, where $\gamma_i \in \mathbb{R}^{N_i}$, such that $\nabla_i L(\mathbf{x}^*, \mu^*) = \nabla_i \Lambda_\rho (\mathbf{x}^*, \mu^*) = \gamma_i, i \in \mathbb{N}_n$. Specifically, when $X = \mathbb{R}^N$, $\gamma = 0$.

\textbf{\textit{Proof:}} See Appendix \ref{subsection: Appendix B}.

From Appendix~\ref{subsection:appendix A}, we have $\mathcal{C}_i^{k} = \phi_i (\bar{\mathbf{x}}_i^{k,v}, \bar{\mathbf{x}}_{-i}^{k,v})$ and $\mathcal{D}_i^{k} = \bar{\mathbf{x}}_i^{k,v} - \bar{\mathbf{x}}_i^{k,v-1}$, which capture the residual information of the solution at the end of the $k$-th outer loop. These definitions apply analogously to the final solution at the end of the outer loop iteration, i.e., the limit point $\mathbf{x}^*$ of $\{\mathbf{x}^k\}$, leading to a similar definition for $\mathcal{C}_i^*$ and $\mathcal{C}^*$.

\noindent \textbf{Theorem 2} Assume \hyperref[assumption:A1]{(A1)}--\hyperref[assumption:A4]{(A4)}. For $k = 0, 1, \ldots$, let $\{\mathbf{x}^{k}\}$ satisfy
\[
    \|\mathcal{D}^{k} \|_\infty \leq \epsilon_{\text{dual}}^k,
\]
and assume that $\{\mu^k\}$ is bounded, and $\{\epsilon_{\text{dual}}^k\}$ and $\{\rho^k\}$ satisfy
\begin{align*}
    0 < \rho_i^k < \rho_i^{k+1}, & \forall k,  \rho_i^k \to \infty, i \in \mathbb{N}_n,    \\
    0 \leq \epsilon_{\text{dual}}^k, & \forall k,  \epsilon_{\text{dual}}^k \to 0.
\end{align*}
Then
\[
    \left\{\mu_i^{k} + 2\rho_i^{k} \circ \rho_i^{k} \circ \mathcal{C}_i^{k}\right\}_{K} \to \mu_i^*,  i \in \mathbb{N}_n,
\]
where $\mu_i^*$ is a vector satisfying, together with $\mathbf{x}^*$, the following conditions
\[
    \nabla L (\mathbf{x}^*, \mu^*)^{\top} (\mathbf{x} - \mathbf{x}^*) \geq 0,  \forall \mathbf{x} \in X, \quad  \mathcal{C}^* = 0.
\]
This means that we obtain the minimum $\mathbf{x}^*$.

\textbf{\textit{Proof:}} Without loss of generality, we assume that the entire sequence $\{\mathbf{x}^k\}$ converges to $\mathbf{x}^*$. Define for all $k$
\begin{equation*}
    \mu_i^{k+1} = \mu_i^k + 2\rho_i^k \circ \rho_i^k \circ \mathcal{C}_i^k,  i \in \mathbb{N}_n.
\end{equation*}
We have
\begin{equation}
    \nabla_i \Lambda_{\rho^k} (\mathbf{x}^k, \mu^k) = \nabla_i f(\mathbf{x}^k) + \left( \nabla_i \mathcal{C}_i^k \right)^{\top} \mu_i^{k+1}, i \in \mathbb{N}_n. \label{eq:17}
\end{equation}
Since we assume that the constraint $\mathcal{C}_i = \phi_i (\mathbf{x}_i, \mathbf{x}_{-i})$ consists of $m_i$ linearly independent constraints, it follows that $\nabla_i \mathcal{C}_i^k$ has rank $m_i$ for all sufficiently large $k$. Then, by multiplying (\ref{eq:17}) with
\begin{equation*}
    \left[\nabla_i \mathcal{C}_i^k\left(\nabla_i \mathcal{C}_i^k\right)^{\top} \right]^{-1} \nabla_i \mathcal{C}_i^k,
\end{equation*}
we obtain
\begin{equation}
\begin{aligned}
\mu_i^{k+1} = & {}
     \left[\nabla_i \mathcal{C}_i^k\left(\nabla_i \mathcal{C}_i^k\right)^{\top} \right]^{-1} \nabla_i \mathcal{C}_i^k [\nabla_i \Lambda_{\rho^k} (\mathbf{x}^k, \mu^k) \\
     & {} - \nabla_i f(\mathbf{x}^k)],  i \in \mathbb{N}_n. \label{eq:18}
\end{aligned}
\end{equation}
By  Theorem 1, it can be concluded that as $k \to \infty$, $\epsilon_{\text{dual}}^k \to \epsilon_{\text{dual}} \approx 0$, there exists a limit point $\mathbf{x}^*$ such that
\begin{equation*}
    \nabla_i \Lambda_{\rho^k} (\mathbf{x}^k, \mu^k) \to \gamma_i,
\end{equation*}
where $\gamma_i \in \mathbb{R}^{N_i},  i \in \mathbb{N}_n$.
Furthermore, (\ref{eq:18}) yields
\begin{equation*}
    \mu_i^{k+1} \to \mu_i^*,  i \in \mathbb{N}_n,
\end{equation*}
where
\begin{equation*}
    \mu_i^* = \left[\nabla_i \mathcal{C}_i^*\left(\nabla_i \mathcal{C}_i^*\right)^{\top} \right]^{-1} \nabla_i \mathcal{C}_i^* \left[ \gamma_i - \nabla_i f(\mathbf{x}^*)\right],i \in \mathbb{N}_n.
\end{equation*}
Using again the fact $\nabla_i \Lambda_{\rho^k} (\mathbf{x}^k, \mu^k) \to  \gamma_i$, along with Lemma 1 and (\ref{eq:17}), we see that
\begin{equation*}
    \left[\nabla_i f(\mathbf{x}^*) + \left(\nabla_i \mathcal{C}_i^*\right)^{\top}  \mu_i^*\right]^{\top} (\mathbf{x}_i - \mathbf{x}_i^*) \geq 0,  \forall \mathbf{x}_i \in X_i,  i \in \mathbb{N}_n.
\end{equation*}
Since $\left\{\mu_i^k,i \in \mathbb{N}_n\right\}$ is bounded and $\mu_i^k + 2\rho_i^k \circ \rho_i^k \circ \mathcal{C}_i^k \to \mu_i^*$, it follows that $\{2\rho_i^k \circ \rho_i^k \circ \mathcal{C}_i^k\}$ is bounded. Since $\rho_i^k \to \infty$, we must have $\mathcal{C}_i^k \to 0$ and we conclude that $\mathcal{C}_i^* = 0$, $i \in \mathbb{N}_n$. 
It implies that $\mathbf{x}^*$ is an optimal solution. \hfill $\square$

Next, we introduce two termination conditions: one equivalent to \hyperref[condition:B2]{(B2)}, and the other more lenient. These conditions facilitate practical implementation and significantly enhance the algorithm's execution efficiency.

\begin{itemize}
    \label{condition:B3}
    \item[(B3)] If the condition $v = v_{\max}^k$ or $\|\mathcal{D}^{k,v}\|_\infty \leq \epsilon_{\text{dual}}$ is satisfied, where $\{v_{\max}^k\}$ satisfies $1 \leq v_{\max}^k$, $\forall k$, $v_{\max}^k \to \infty$, proceed to Step 2.
\end{itemize}
\begin{itemize}
    \label{condition:B4}
    \item[(B4)] If the condition $v = v_{\max}$ or $\|\mathcal{D}^{k,v}\|_\infty \leq \epsilon_{\text{dual}}$ is satisfied, proceed to Step 2.
\end{itemize}

For condition \hyperref[condition:B4]{(B4)}, it is clearly understood that it is not equivalent to \hyperref[condition:B2]{(B2)}. Furthermore, when setting the termination conditions \hyperref[condition:B2]{(B2)}, \hyperref[condition:B3]{(B3)}, and \hyperref[condition:B4]{(B4)},  it is crucial to update the stopping criteria for the outer loop as follows:
\begin{equation*}
     \text{If } \|\mathcal{C}^{k,v}\|_\infty \leq \epsilon_{\text{pri}} \text{ and } \|\mathcal{D}^{k,v}\|_\infty \leq \epsilon_{\text{dual}} \text{ are satisfied},
\end{equation*}
to prevent the algorithm from terminating prematurely before $\mathcal{D}^{k,v}$ meets the specified conditions. Here, $\|\mathcal{C}^{k,v}\|_\infty=\max_{i \in \mathbb{N}_n} \| \mathcal{C}_i^{k,v} \|_\infty$, and $\| \mathcal{C}_i^{k,v} \|_\infty = \max_{j \in \mathbb{N}_{m_i}} |\mathcal{C}_{ij}^{k,v}|$.

It is worth noting that all the parameter settings mentioned earlier can be freely adjusted according to specific requirements or experimental observations.

\begin{table*}[t!htb]
    \centering
    \begin{minipage}{0.95\textwidth} 
        \centering
                \caption{Comparison of the Distributed Algorithms Discussed}
        \label{tab:1}
        \renewcommand{\arraystretch}{1.6}
        \setlength{\tabcolsep}{1.5pt} 
        \footnotesize 
        \begin{tabular*}{\linewidth}{@{\extracolsep{\fill}} cllccccc }
            \toprule
            Problem types & Methods & Problem structures & Main assumptions & Termination criteria & $h$-level & References \\
            \midrule
            \multirow{3.5}{*}{unconstrained}             
            & BCD 
            & \tabincell{l}{$\min ~ f(\mathbf{x}_1, \mathbf{x}_2, \dots, \mathbf{x}_n)$ 
            } 
            & \tabincell{c}{ $f$ continuously  \\ differentiable  convex }
            & \hyperref[condition:B1]{(B1)}  & $\geq 2$ & \cite{bertsekasConvexOptimizationAlgorithms2015} \\

            \cmidrule{2-7}
            
            & BSUM 
            & \tabincell{l}{$\min ~ \sum_{i=1}^n f_i(\mathbf{x}_i) + f_0(\mathbf{x})$ 
            } 
            & \tabincell{c}{$f_0$ smooth convex, \\ $f_i$ nonsmooth convex }
            & \hyperref[condition:B1]{(B1)} & $\geq 2$ & \cite{hongIterationComplexityAnalysis2017,razaviyayn2013unified} \\ 

            \cmidrule{1-7}
            
            \multirow{14.5}{*}{constrained} 
            & DQA 
            & \tabincell{l}{$\min ~ \sum_{i=1}^n f_i(\mathbf{x}_i)$ \\~s.t.  $\sum_{i=1}^n \mathbf{A}_i \mathbf{x}_i = \mathbf{b}$, \\ ~~~~ $\mathbf{x}_i \in X_i, i \in \mathbb{N}_n$.}
            & $f_i$ convex 
            & \hyperref[condition:B1]{(B1)} & $2$ & \cite{ruszczynski1995convergence} \\

            \cmidrule{2-7}
            
            & ADAL 
            & \tabincell{l}{$\min ~ \sum_{i=1}^n f_i(\mathbf{x}_i)$ \\~s.t.  $\sum_{i=1}^n \mathbf{A}_i \mathbf{x}_i = \mathbf{b}$, \\ ~~~~ $\mathbf{x}_i \in X_i, i \in \mathbb{N}_n$.} 
            & $f_i$ convex 
            & \hyperref[condition:B4]{(B4)} with  $v_{\text{max}} = 1$ & $2$ & \cite{chatzipanagiotis2015augmented} \\

            \cmidrule{2-7}
            
            &ADMM 
            & \tabincell{l}{$\min ~ \sum_{i=1}^2 f_i(\mathbf{x}_i)$ \\~s.t.  $\sum_{i=1}^2 \mathbf{A}_i \mathbf{x}_i = \mathbf{b}$, \\ ~~~~ $\mathbf{x}_i \in X_i, i \in \mathbb{N}_2$.}
            & $f_i$ convex 
            & \hyperref[condition:B4]{(B4)} with  $v_{\text{max}} = 1$ & $2$ & \cite{boyd2011distributed} \\

            \cmidrule{2-7}
            
            & BSUM-M 
            & \tabincell{l}{$\min ~ \sum_{i=1}^n f_i(\mathbf{x}_i) + f_0(\mathbf{x})$ \\~s.t. $\sum_{i=1}^n \mathbf{A}_i \mathbf{x}_i = \mathbf{b}$, \\~~~~ $\mathbf{x}_i \in X_i, i \in \mathbb{N}_n$.}
            & \tabincell{c}{$f_0$ smooth convex, \\ $f_i$ nonsmooth convex } 
            & \hyperref[condition:B4]{(B4)} with  $v_{\text{max}} = 1$ & $\geq 2$ & \cite{hong2020block} \\

            \cmidrule{2-7}
            
            & DALD 
            & \tabincell{l}{ $\min ~ f(\mathbf{x}_1, \mathbf{x}_2, \dots, \mathbf{x}_n)$ \\ ~s.t. $\sum_{i=1}^n \mathbf{A}_i \mathbf{x}_i = \mathbf{b}$, \\ ~~~~~$\mathbf{x}_i \in X_i, i \in \mathbb{N}_n$.} 
            & \tabincell{c}{ $f$  continuously  \\ differentiable  convex\footnotemark  }
            & \tabincell{c}{\hyperref[condition:B1]{(B1)}, \hyperref[condition:B2]{(B2)}, \\ \hyperref[condition:B3]{(B3)}, or \hyperref[condition:B4]{(B4)}} & $\geq 2$ & This paper \\
            \bottomrule
        \end{tabular*}
    \end{minipage}
\end{table*}

\subsection{Framework Unification}

A wide variety of distributed constrained optimization algorithms based on the AL techniques currently exist. This section aims to provide a unified analysis of these methods and explore their common characteristics. However, before undertaking this integration, we first need to re-examine the principles and mechanisms of DALD. Specifically, when addressing a constrained optimization problem (\ref{eq:1}) (we will temporarily disregard the convex set constraint $\mathbf{x} \in X$ in this section), the process typically involves the following three steps:

\begin{itemize}
    \item[1)] \textit{Decomposition}: Identify local variables and coupling relationships according to decision requirements, and decompose the original problem into multiple subproblems.
    \item[2)] \textit{Relaxation}: Use the AL technique to convert the subproblems into unconstrained forms of the local augmented Lagrangian.
    \item[3)] \textit{Coordination}: Determine the subproblem solving sequence, and update variables and related parameters during the iterative process.
\end{itemize}

\begin{figure}[t!]
\centering
 \includegraphics[width=1\linewidth]{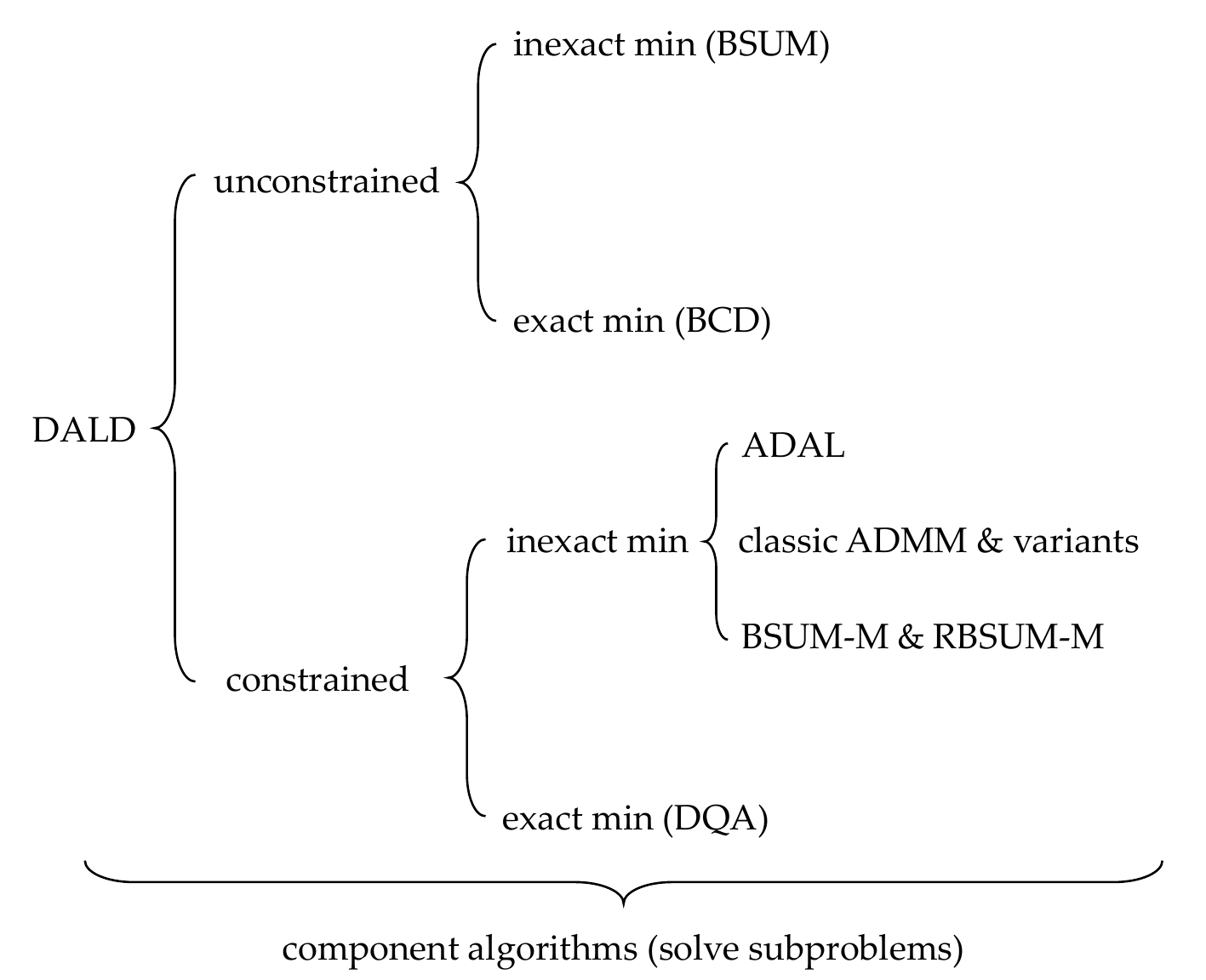}
\caption{A unified framework for distributed optimization}
\label{fig:4}
\end{figure}

The relaxation process introduces Lagrange multipliers, which are updated by the outer layer loop. During the coordination process, the local variables are alternately optimized in the inner layer loop, implying that the number of levels $h > 1$ in the hierarchical coordination network. When the constraint set $\varphi(\mathbf{x})$ is empty, the standard version of DALD reduces to the well-known BCD. Thus, from the perspective of decomposition and coordination, BCD can also be viewed as a form of distributed optimization framework.

For typical BCD methods, different component algorithms can be selected in the solver layer based on the characteristics of the subproblems. 
For instance, if the subproblems are differentiable, gradient-based algorithms such as classic Gradient Descent (GD) or Newton's method (NM) can be considered. For example, in the case of the classic GD (see Appendix \ref{subsection: Appendix C} for details), the gradient of the objective function with respect to the elements $x_{ij}$ in block $\mathbf{x}_i$ needs to be computed, involving all information related to $x_{ij}$, $i \in \mathbb{N}_1$, $j \in \mathbb{N}_{N_i}$. From the perspective of decomposition and coordination, gradient descent algorithms can also be viewed as a generalized distributed optimization algorithm \cite{cohen1978optimization}. According to (\ref{eq:VII-C-3}), first-order optimization algorithms that use gradient information (such as GD) are implemented in the solver layer in a manner consistent with Jacobian-type standards.

In designing accelerated versions of DALD, although an inexact approach is used in the inner layer to solve the augmented Lagrangian function, in the solver layer we generally assume that the subproblems can be solved exactly or disregard this detail. Most existing distributed optimization algorithms based on augmented Lagrangian techniques follow this principle, including ADAL, classic ADMM and its various variants. Since these algorithms set the maximum number of iterations in the inner layer to $v_{\max}=1$, they may not guarantee convergence when dealing with more general optimization problems \cite{chen2016direct}.

Based on the previous discussion, a natural question arises: can the DALD still maintain convergence when subproblems are solved inexactly in the solver layer? The answer is affirmative. A method called Block Successive Upper-bound Minimization (BSUM), introduced in \cite{hongIterationComplexityAnalysis2017,razaviyayn2013unified}, serves as an inexact version of the BCD method and ensures convergence.  Similarly, if both the inner and solver layers are solved in an inexact manner, convergence can still be guaranteed as long as the precision of the solution is gradually increased during iterations \cite[Ch. 3.4.4]{boyd2011distributed}. Block Successive Upper-bound Minimization Method of Multipliers (BSUM-M) and Randomized  Block Successive Upper-bound Minimization Method of Multipliers (RBSUM-M) are based on this principle\cite{hong2020block}, with the distinction that the former employs full-cycle coordination in the inner layer, while the latter uses partial-cycle coordination.

As illustrated in Fig.~\ref{fig:4}, we propose a unified distributed optimization framework based on the AL techniques---DALD. Unlike the traditional focus on only the outer and inner layers, we introduce the concept of the ``solver layer,'' which involves component algorithms   responsible for iterative solving of the subproblems. From this perspective, first-order algorithms and ALM (when first-order algorithms are used to solve (\ref{eq:4})) can be viewed as generalized distributed optimization methods. These methods  decompose $f(\mathbf{x})$  into at least \(j \geq 2\) subproblems within the solver layer and solve them through coordination, but essentially, coordination still occurs within a single system, i.e., $n=1$. However, unless otherwise specified, the term ``distributed optimization" in this paper typically refers to optimization frameworks in the context of constrained optimization problems, characterized by the requirement that the entire system must involve at least  \(n \geq 2\) agents.

For clarity, Table~\ref{tab:1} provides a comparison to highlight the differences among the involved algorithms. It can be observed that DALD handles a more general range of problem forms and is more versatile in terms of its termination criteria. 
It is important to note that the algorithms shown in Fig.~\ref{fig:4} are mature classical algorithms proposed by previous researchers, and their variants are not depicted in the figure. We believe that the problem forms and coordination mechanisms addressed by these variants are essentially similar to those of the prototype algorithms and fall within the interpretive scope of DALD. For instance, most variants of the ADMM are still confined to problems with fully separable and convex objectives,  or 2-level coordination. 

\footnotetext{The convergence of the non-differentiable version has been proven in our latest work~\cite{guo2025distributedoptimizationdesignedfederated}.}

As a distributed optimization framework with three layers of loops, DALD is not limited to the algorithms listed in Fig.~\ref{fig:4} but can be flexibly adjusted and customized according to different optimization needs, thereby fully leveraging its powerful adaptability.

\section{Numerical Experiments}
\label{sec:5.Numerical Experiments}
This section presents the numerical results of applying the DALD method to the optimization problems, aiming to further illustrate the proposed approach. For convenience, we set the initial parameters as $\mu^1 = \mathbf{0}_m$, $\rho = \mathbf{1}_m$, and the initial solution as $\mathrm{x}^{1,0} = \mathbf{0}_N$. Unless specified otherwise, we assume $\epsilon_{\text{dual}} = \epsilon_{\text{pri}} = 1 \times 10^{-3}$ by default.

In the solver layer, we assume exact solutions are obtained whenever a third-party solver is invoked and successfully returns a solution. Furthermore, our algorithm is implemented in Python 3.12.4, and all experiments were conducted on a 64-bit Windows PC equipped with an Intel i5-13500H processor (2.60 GHz) and 32.0 GB of RAM.

\begin{figure*}[t!bp]
\centering
\includegraphics[scale=0.33]{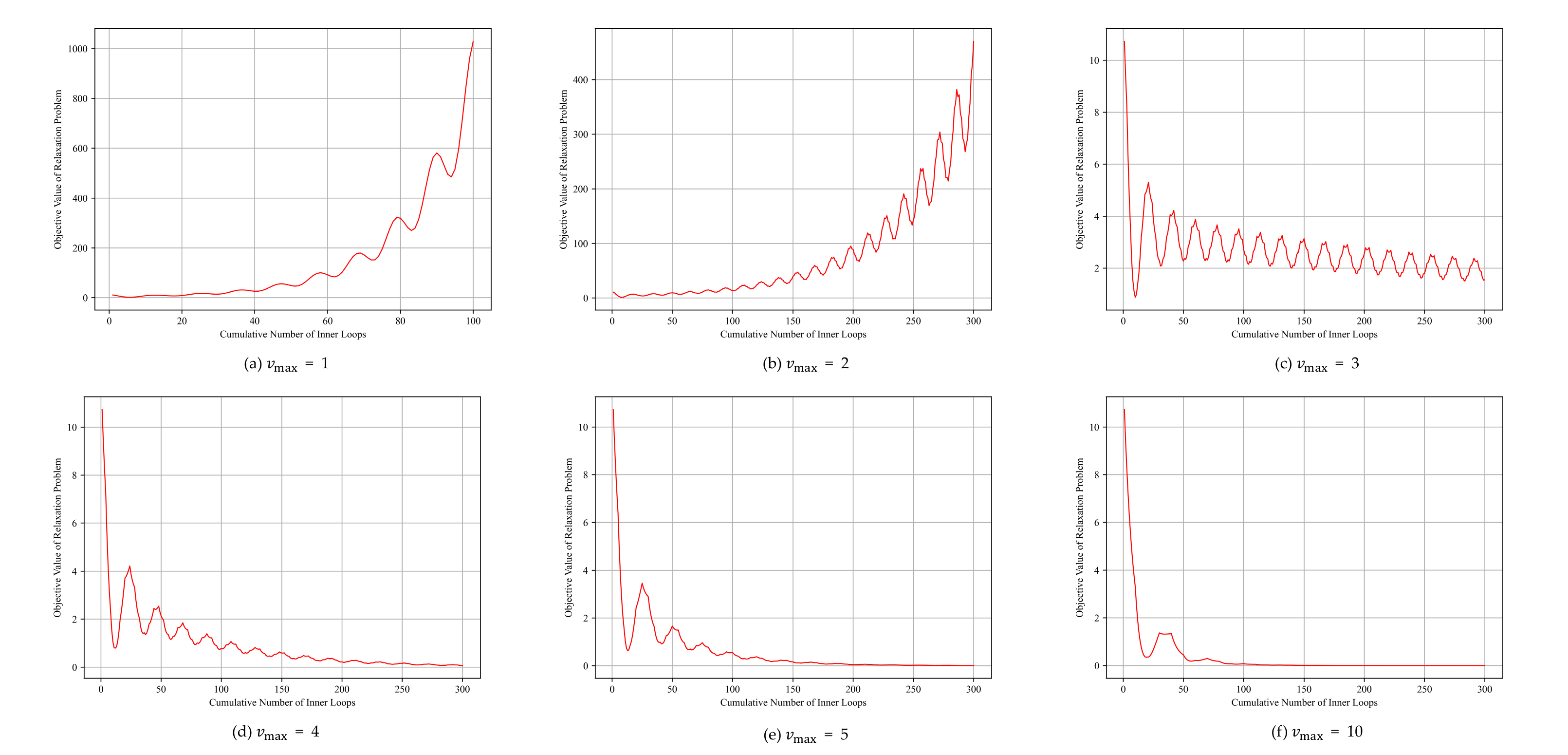}
\caption{Convergence trends for different values of $v_{\max}$}
\label{fig:5}
\end{figure*}

\subsection{A Case of ADMM Non-Convergence}
Classic ADMM does not necessarily guarantee convergence when using its direct extension for sequentially alternating optimization of three or more blocks of variables, that is, the number of hierarchical coordination levels is \( h \geq 3 \) in the inner layer. The reasons for the potential lack of convergence are readily apparent from Section \ref{subsection:4.2 Parameters Determination and Accelerated Version}. To further validate our findings, we analyze a counterexample where an extended version of ADMM fails to converge \cite{chen2016direct}. Consider the following optimization problem:

\begin{equation*}
\min ~0
\end{equation*}
\vspace{-1.5em}
\begin{equation*}
\text{s.t.~} A_1 x_1 + A_2 x_2 + A_3 x_3 = 0,
\end{equation*}
where \( A = \begin{bmatrix}
1 & 1 & 1 \\
1 & 1 & 2 \\
1 & 2 & 2
\end{bmatrix} \). For this problem, we can readily obtain the AL function and subsequently perform sequential minimization over \( x_1 \), \( x_2 \), and \( x_3 \). By applying optimality conditions in the solver layer, we obtain analytical solutions for the inner layer
\begin{equation*}
x_1^{k,v+1} = -\frac{A_1^{\top}}{\|A_1\|^2}  \left( \frac{\mu^k}{2\rho \circ \rho} + A_2 x_2^{k,v} + A_3 x_3^{k,v} \right),
\end{equation*}
\begin{equation*}
x_2^{k,v+1} = -\frac{A_2^{\top}}{\|A_2\|^2}  \left( \frac{\mu^k}{2\rho \circ \rho} + A_1 x_1^{k,v+1} + A_3 x_3^{k,v} \right),
\end{equation*}
\begin{equation*}
x_3^{k,v+1} = -\frac{A_3^{\top}}{\|A_3\|^2}  \left( \frac{\mu^k}{2\rho \circ \rho} + A_1 x_1^{k,v+1} + A_2 x_2^{k,v+1} \right).
\end{equation*}

If the stopping condition for the DALD inner loop uses \hyperref[condition:B4]{(B4)} and \( v_{\max} = 1 \), we obtain the direct extension of classic ADMM. As shown in Fig.~\ref{fig:5}(a), the algorithm does not converge under this condition. Consequently, we gradually increased the value of \( v_{\max} \). As illustrated in Fig.~\ref{fig:5}, when \( v_{\max} \geq 3 \), the algorithm begins to exhibit convergence behavior. Although only the first 300 cumulative inner loop iterations are shown in the figure for comparative purposes (except for \( v_{\max} = 1 \)), results ultimately converge to the optimal solution \( x = [0, 0, 0] \) when \( v_{\max} \geq 3 \), and the AL function value also converges to 0.

\begin{figure}[tbp]
\centering
\includegraphics[scale=0.32]{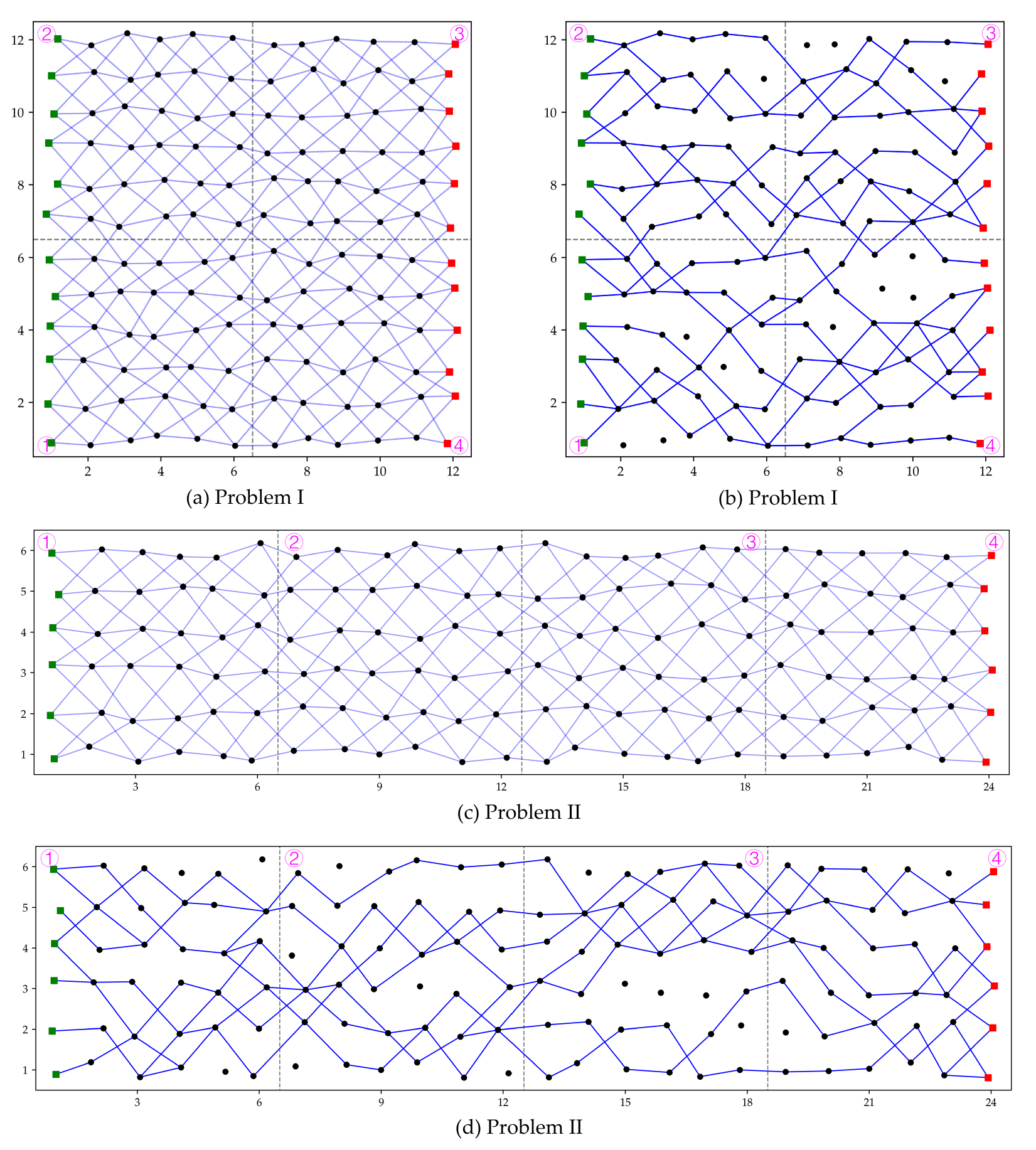}
\caption{Randomly generated networks for the LNF problem}
\label{fig:7}
\end{figure}

\subsection{Network Flow Problem}
We now turn our attention to the Linear Network Flow (LNF) problem, which involves minimizing costs associated with linear arc costs. This problem is central to many optimization scenarios, arising naturally in the design and analysis of large-scale systems, including supply chains, energy distribution, communication, transportation, and manufacturing networks~\cite{bertsekasNetworkOptimizationContinuous1998}.

Consider a directed graph $G = (\mathcal{N}, A)$, where $\mathcal{N}$ is the set of vertices and $A$ is the set of directed edges.  Each edge \( (i,j) \) has an associated cost \( c_{ij} \), and the flow through this edge is denoted by \( t_{ij} \). The flow \( t_{ij} \) must satisfy constraints \( a_{ij} \leq t_{ij} \leq b_{ij} \). The supply at node \( i \) is given by \( s_i \), where \( a_{ij} \), \( b_{ij} \), and \( s_i \) are predefined. For source nodes \( S \), we have \( s_i > 0 \) for \( \forall i \in S \), while for sink nodes \( D \), \( s_i < 0 \), \( \forall i \in D \). Consider the problem:

\begin{equation*}
\min \sum_{(i,j) \in A} c_{ij} t_{ij}
\end{equation*}
\begin{equation*}
\text{s.t.~} \sum_{\{j | (i,j) \in A\}} t_{ij} - \sum_{\{j | (j,i) \in A\}} t_{ji} = s_i,  \forall i \in \mathcal{N}
\end{equation*}
\begin{equation*}
a_{ij} \leq t_{ij} \leq b_{ij},  \forall (i,j) \in A
\end{equation*}

In the specific instances considered, we define a subset \( T \) of nodes as transshipment nodes, where \( s_i = 0 \), \( \forall i \in T \). Additionally, we model \( s_i \) for nodes in \( S \cup D \) using a Poisson distribution with a mean of 50. The flow bounds on each arc are set between 0 and 50, and the cost coefficients are assigned random integer values between 1 and 5. In Fig.~\ref{fig:7}(a) and \ref{fig:7}(c), two randomly generated LNF problems (Problem I and Problem II) with a 144-node configuration are presented. Green dots represent source nodes, red dots indicate sink nodes, and black dots denote transshipment nodes. All available arcs are illustrated as blue lines.

\begin{figure*}[t!bp]
\centering
\includegraphics[scale=0.33]{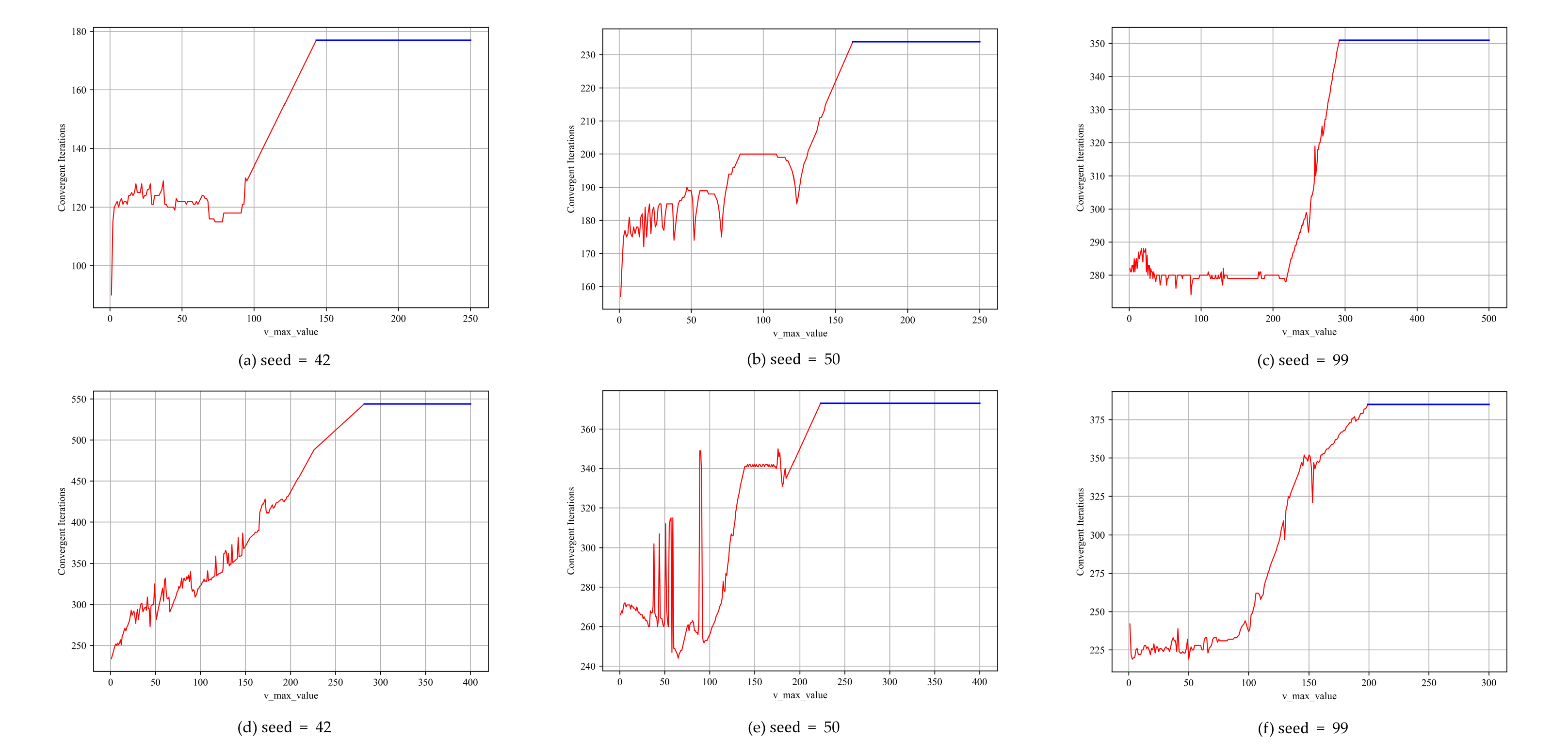}
\caption{Required convergent iterations for diffirent scenarios of the LNF problems}
\label{fig:8}
\end{figure*}

To simplify the process, the LNF problems were divided into four equally-sized, interconnected subproblems, respectively, as depicted in Fig.~\ref{fig:7}(a) and \ref{fig:7}(c). Leveraging the inherent coupling among these subproblems, we employed the Gurobi optimizer (version 11.0.3) to solve the subproblems in the solver layer, while the specific subproblem-solving sequence 
\[
\mathbb{S} = \begin{bmatrix}
3 & 1 & 1 & 1 \\
0 & 2 & 1 & 1 \\
0 & 0 & 1 & 1 \\
0 & 0 & 0 & 0
\end{bmatrix}
\]
was applied in the inner layer. The inner loop termination condition follows \hyperref[condition:B4]{(B4)}, and we record the total number of inner loop iterations required for different \( v_{\max} \) values. The results demonstrate that, irrespective of the \( v_{\max} \) value, the cases consistently converge. Fig.~\ref{fig:8}(a-c) and \ref{fig:8}(d-f) show the iteration counts for the two problems under different random seed settings. The blue line in each figure indicates that the algorithm satisfies the standard version of DALD, when \( v_{\max} \) is set to the corresponding value. The minimum iteration cost generally occurs at smaller \( v_{\max} \) values, though not always at \( v_{\max} = 1 \). By selecting an appropriate \( v_{\max} \), we observe that the accelerated DALD can significantly reduce computational costs compared to the standard version. Additionally, Fig.~\ref{fig:7}(b) and \ref{fig:7}(d) also display the flow routing decisions obtained after solving the problem using DALD with a random seed value of 50.

\section{Conclusion}
\label{sec:6.Conclusion}
This paper investigates a distributed, solver-agnostic optimization technique, DALD, which is distinguished by its ability to solve optimization problems with broader characteristics under more flexible communication topologies. Specifically, we introduce both the standard and accelerated versions of the method and provide the corresponding convergence proofs. In addition, we indirectly address a longstanding issue in the field of distributed optimization: why the direct extension of ADMM is not necessarily convergent. Unlike previous research, we introduce the concept of ``solver layer" and highlight that DALD can leverage different solvers to solve subproblems without depending on any specific solver.  We also briefly review and analyze existing distributed optimization methods, revealing that they fall within the explanatory scope of the DALD framework. Different from Professor Bertsekas' perspective \cite[Ch. 7.4.2]{bertsekas_nonlinear_2016}, we view DALD as a distributed implementation of alternating optimization within the framework of the AL for constrained optimization problems.

\section{Appendices}
\label{sec:7.Appendices}
\subsection{Proof of the Derivation about \( \mathcal{C}^{\textit{k,v}} \) and \( \mathcal{D}^{\textit{k,v}} \)}
\label{subsection:appendix A}

The necessary and sufficient optimality conditions for the problem (\ref{eq:6}) are primal feasibility,
\begin{equation}
\phi_i (\mathbf{x}_i^*, \mathbf{x}_{-i}^*) = 0,  i \in \mathbb{N}_n \label{eq:VII-A-1}
\end{equation}
and dual feasibility deduced by problem (\ref{eq:8}),
\begin{equation}
\nabla_i \mathcal{L}_i (\mathbf{x}_i^*, \mathbf{x}_{-i}^*, \mu_i^*)^{\top} \left(\mathbf{x}_i - \mathbf{x}_i^*\right) \geq 0, \forall \mathbf{x}_i \in X_i, i \in \mathbb{N}_n.
 \label{eq:VII-A-2}
\end{equation}
According to Theorem 1, the optimal point \( \bar{\mathbf{x}}^{k,v} \), as the limit point of the sequence \(\{\mathbf{x}^{k,v}\}\), is also a critical point of \(\Lambda_\rho (\mathbf{x}, \mu^k)\). This ensures that the iterative process reaches a steady state. Therefore, when \(v \to \infty\), resulting in
\begin{equation}
\mathcal{D}_i^{k,v} = \mathbf{x}_i^{k,v} - \mathbf{x}_i^{k,v-1} = 0,   \label{eq:VII-A-3}
\end{equation}
which can be viewed as a residual for the dual feasibility condition (\ref{eq:VII-A-2}) at the \((k,v)\)-th loop iteration, $i \in \mathbb{N}_n$.

Let \(\mathbf{x}_i^{k,v}\) minimize \(\mathcal{A}_{\rho_i}^i (\mathbf{x}_i, \mathbf{w}_{-i}^{k,v}, \mu_i^k)\) by definition, \(i \in \mathbb{N}_n\), we have that
\begin{equation}
\nabla_i \mathcal{A}_{\rho_i}^i ( \mathbf{x}_i^{k,v}, \mathbf{w}_{-i}^{k,v}, \mu_i^k)^{\top} (\mathbf{x}_i - \mathbf{x}_i^{k,v}) \geq 0, \forall \mathbf{x}_i \in X_i,
 \label{eq:VII-A-4}
\end{equation}
where the gradient of the local augmented Lagrangian $\nabla_i \mathcal{A}_{\rho_i}^i ( \mathbf{x}_i^{k,v}, \mathbf{w}_{-i}^{k,v}, \mu_i^k)$ 
\begin{multline*}
= \nabla_i f_i (\mathbf{x}_i^{k,v}, \mathbf{w}_{-i}^{k,v}) + \nabla_i \phi_i (\mathbf{x}_i^{k,v}, \mathbf{w}_{-i}^{k,v})^{\top}\mu^k\\
+ \nabla_i \phi_i (\mathbf{x}_i^{k,v}, \mathbf{w}_{-i}^{k,v})^{\top} \left[2\rho_i \circ \rho_i \circ \phi_i (\mathbf{x}_i^{k,v}, \mathbf{w}_{-i}^{k,v})\right].
\end{multline*}
Given that \(\mathbf{w}_{-i}^{k,v} = (\mathbf{x}_1^{k,v}, \ldots, \mathbf{x}_{i-1}^{k,v}, \mathbf{x}_{i+1}^{k,v-1}, \ldots, \mathbf{x}_n^{k,v-1})\), after the last subproblem is resolved, \(\mathbf{x}\) will be updated as  \(\mathbf{x}^{k,v} = (\mathbf{x}_1^{k,v}, \mathbf{x}_2^{k,v}, \ldots, \mathbf{x}_n^{k,v})\).
When \(v \to \infty\), we obtain the optimal point \( \bar{\mathbf{x}}^{k,v} (\mathbf{x}^{k,*}) \) and $\nabla_i \mathcal{A}_{\rho_i}^i (\bar{\mathbf{x}}_i^{k,v}, \bar{\mathbf{x}}_{-i}^{k,v}, \mu_i^k)$
\begin{equation}
\begin{aligned}
\quad = {} & \nabla_i f_i (\bar{\mathbf{x}}_i^{k,v}, \bar{\mathbf{x}}_{-i}^{k,v}) + \nabla_i \phi_i (\bar{\mathbf{x}}_i^{k,v}, \bar{\mathbf{x}}_{-i}^{k,v})^{\top}\mu^{k} \\ 
& + \nabla_i \phi_i (\bar{\mathbf{x}}_i^{k,v}, \bar{\mathbf{x}}_{-i}^{k,v})^{\top} \left[2\rho_i \circ \rho_i \circ \phi_i (\bar{\mathbf{x}}_i^{k,v}, \bar{\mathbf{x}}_{-i}^{k,v})\right] \\
= {} & \nabla_i f_i (\bar{\mathbf{x}}_i^{k,v}, \bar{\mathbf{x}}_{-i}^{k,v}) + \nabla_i \phi_i (\bar{\mathbf{x}}_i^{k,v}, \bar{\mathbf{x}}_{-i}^{k,v})^{\top}\mu^{k+1}. 
\label{eq:VII-A-5}
\end{aligned}
\end{equation}
As \(k \to \infty\), we have \( \bar{\mathbf{x}}^{k,v}(\bar{\mathbf{x}}^{k}) \to \mathbf{x}^* \) and \( \bar{\mu}^k \to \mu^*\), which means that if the \(k\)-th outer loop of the DALD converges, i.e., satisfying condition (\ref{eq:VII-A-3}), then \( \bar{\mathbf{x}}_i^{k,v} \) and \( \mu_i^{k+1} \) always satisfy (\ref{eq:VII-A-2}), \(i \in \mathbb{N}_n\).
We will refer to
\begin{equation}
\mathcal{C}_i^k = \mathcal{C}_i^{k,v} = \phi_i (\bar{\mathbf{x}}_i^{k,v}, \bar{\mathbf{x}}_{-i}^{k,v}) \label{eq:VII-A-6}
\end{equation}
as the \textit{primal residual} and 
\begin{equation}
\mathcal{D}_i^k = \mathcal{D}_i^{k,v} = \bar{\mathbf{x}}_i^{k,v} - \bar{\mathbf{x}}_i^{k,v-1} \label{eq:VII-A-7}
\end{equation}
as the \textit{dual residual} at outer loop \(k\) in DALD, \(i \in \mathbb{N}_n\).

\subsection{Proof of Lemma 1}
\label{subsection: Appendix B}
\textbf{\textit{Proof:}}  Without loss of generality, we reformulate the constraint as $X = \{\mathbf{x} \in \mathbb{R}^N : h(\mathbf{x}) = 0\}$ of (\ref{eq:1}), where $h : \mathbb{R}^{N} \to \mathbb{R}^p$, and let $\lambda \in \mathbb{R}^p$ be the Lagrange multipliers associated with $h(\mathbf{x})$. According to the KKT conditions, we have
\begin{equation*}
\varphi(\mathbf{x}^*) = 0, \quad h(\mathbf{x}^*) = 0,
\end{equation*}
and
\begin{equation*}
0 = \nabla_i f(\mathbf{x}^*) + \nabla_i \varphi(\mathbf{x}^*)^{\top} \mu^* + \nabla_i h(\mathbf{x}^*)^{\top} \lambda^*,  i \in \mathbb{N}_n.
\end{equation*}
Rearranging gives:
\begin{equation*}
-\nabla_i h(\mathbf{x}^*)^{\top} \lambda^* = \nabla_i f(\mathbf{x}^*) + \nabla_i \varphi(\mathbf{x}^*)^{\top} \mu^*,  i \in \mathbb{N}_n.
\end{equation*}
Letting $\gamma_i = -\nabla_i h(\mathbf{x}^*)^{\top} \lambda^*$, we obtain $\nabla_i L(\mathbf{x}^*, \mu^*) = \gamma_i$. Since $\varphi(\mathbf{x}^*)=0$, it follows that $\nabla_i L(\mathbf{x}^*, \mu^*) = \nabla_i \Lambda_\rho (\mathbf{x}^*, \mu^*)$. When $X = \mathbb{R}^N$, the conclusion holds naturally according to the KKT conditions. \hfill $\square$

\subsection{Gradient Descent and BCD}
\label{subsection: Appendix C}

To facilitate the exposition using the GD algorithm, we assume that all constraints for the problem (\ref{eq:1}) are absent, forming an empty set, and that the function \( f \) is differentiable. We obtain an unconstrained optimization problem:
\begin{equation}
\min f(\mathbf{x}_1, \mathbf{x}_2, \ldots, \mathbf{x}_n) \label{eq:VII-C-1}
\end{equation}




Given the current iterate \( \mathbf{x}^r = (\mathbf{x}_1^r, \mathbf{x}_2^r, \ldots, \mathbf{x}_n^r) \), we generate the next iterate \( \mathbf{x}^{r+1} = (\mathbf{x}_1^{r+1}, \mathbf{x}_2^{r+1}, \ldots, \mathbf{x}_n^{r+1}) \). 
When directly using GD to optimize (\ref{eq:VII-C-1}), the iteration formula in the solver layer is as follows for each $i \in \mathbb{N}_1$:
\begin{equation}
\begin{aligned}
\mathbf{x}_i^{r+1} = \mathbf{x}_i^r - \alpha_r \nabla_i f(\mathbf{x}^r)  = \mathbf{x}_i^r - \alpha_r \nabla_i f_i (\mathbf{x}_i^r, \mathbf{x}_{-i}^r). 
\end{aligned}
\label{eq:VII-C-2}
\end{equation}

For specific elements,
\begin{equation}
x_{ij}^{r+1} = x_{ij}^r - \alpha_r \nabla_{ij} f_i (\mathbf{x}_i^r, \mathbf{x}_{-i}^r),  i \in \mathbb{N}_1,  j \in \mathbb{N}_{N_i}. \label{eq:VII-C-3}
\end{equation}

The above update rule is equivalent to solving the subproblem:
\begin{equation}
\begin{aligned}
\mathbf{x}_i^{r+1} = {} & \arg\min_{\mathbf{x}_i}  f(\mathbf{x}^r) + \nabla_i f(\mathbf{x}^r) (\mathbf{x}_i - \mathbf{x}_i^r) + \frac{1}{2 \alpha_r} (\mathbf{x}_i - \mathbf{x}_i^r)^2  \\ 
= {} & \arg\min_{\mathbf{x}_i}  \nabla_i f_i (\mathbf{x}_i^r, \mathbf{x}_{-i}^r) (\mathbf{x}_i - \mathbf{x}_i^r)    + \frac{1}{2 \alpha_r} (\mathbf{x}_i - \mathbf{x}_i^r)^2 . 
\end{aligned}
\label{eq:VII-C-4}
\end{equation}


Given the current iterate \( \mathbf{x}^v = (\mathbf{x}_1^v, \mathbf{x}_2^v, \ldots, \mathbf{x}_n^v) \), we generate the next iterate \( \mathbf{x}^{v+1} = (\mathbf{x}_1^{v+1}, \mathbf{x}_2^{v+1}, \ldots, \mathbf{x}_n^{v+1}) \). The typical implementation of the BCD method proceeds as follows:
\begin{equation}
\begin{aligned}
\mathbf{x}_i^{v+1} = {} & \arg\min_{\mathbf{x}_i} f(\mathbf{x}_1^{v+1}, \ldots, \mathbf{x}_{i-1}^{v+1}, \mathbf{x}_i, \mathbf{x}_{i+1}^v, \ldots, \mathbf{x}_n^v) \\ 
= {} & \arg\min_{\mathbf{x}_i} f_i (\mathbf{x}_i, \mathbf{w}_{-i}^v),  i \in \mathbb{N}_n. \label{eq:VII-C-5}
\end{aligned}
\end{equation}

In the solver layer of BCD, with \( \mathbf{w}_{-i}^v \) fixed, we proceed to employ GD iteratively to solve subproblem (\ref{eq:VII-C-5}) until convergence, thereby obtaining \( \mathbf{x}_i^{v+1} \):
\begin{equation}
\mathbf{x}_i^{v,r+1} = \mathbf{x}_i^{v,r} - \alpha_r \nabla_i f_i (\mathbf{x}_i^{v,r}, \mathbf{w}_{-i}^v), i \in \mathbb{N}_n. \label{eq:VII-C-6}
\end{equation}

In the solver layer of DALD, with \( \mathbf{w}_{-i}^{k,v} \) fixed, we employ GD iteratively to solve subproblem (\ref{eq:12}) until convergence, thereby obtaining \( \mathbf{x}_i^{k,v+1} \):
\begin{equation}
\mathbf{x}_i^{k,v,r+1} = \mathbf{x}_i^{k,v,r} - \alpha_r \nabla_i \mathcal{A}_{\rho_i}^i (\mathbf{x}_i^{k,v,r}, \mathbf{w}_{-i}^{k,v}, \mu_i^k),  i \in \mathbb{N}_n. \label{eq:VII-C-7}
\end{equation}

When there are convex set constraints, applicable methods still exist, such as combining with projection operations.




\bibliography{references.bib} 

\bibliographystyle{IEEEtran}

\begin{IEEEbiography}[{\includegraphics[width=1in,height=1.25in,clip,keepaspectratio]{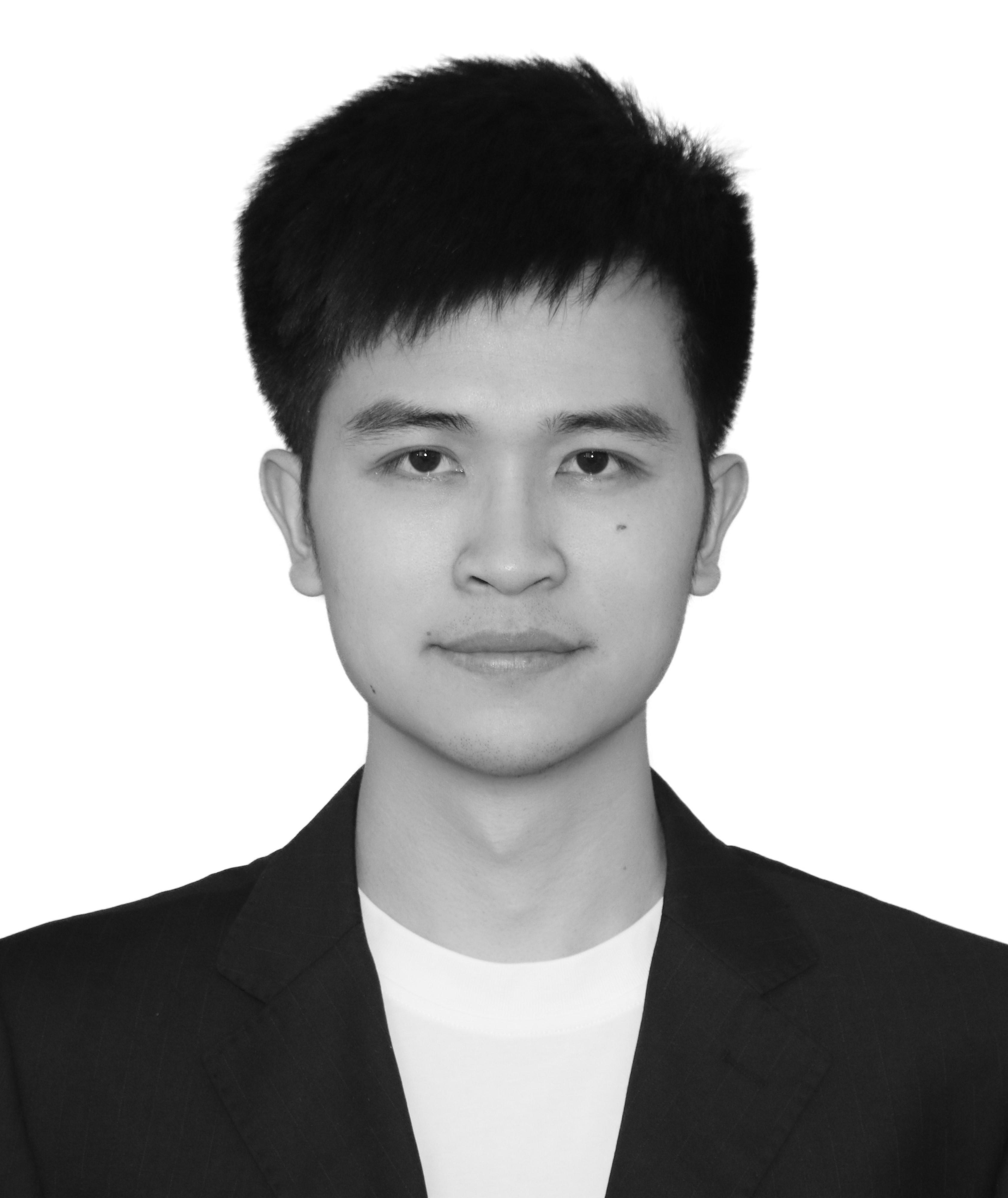}}]{Wenyou Guo}  received his B.Eng. and M.Eng. degrees  in Industrial Engineering from Jiangxi University of Science and Technology, Ganzhou, China, in 2019 and 2022, respectively. He is currently pursuing a Ph.D. degree in Management Science and Engineering at Jinan University, Guangzhou, China. His current research interests include distributed optimization,  federated learning, blockchain, and intelligent manufacturing.
\end{IEEEbiography}



\begin{IEEEbiography}
[{\includegraphics[width=1in,height=1.25in,clip,keepaspectratio]{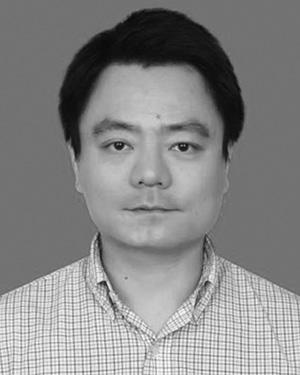}}]{Ting Qu}  received the B.Eng. and M.Phil. degrees in Mechanical Engineering from Xi’an Jiaotong University, Xi’an, China, in 2001 and 2004, respectively, and the Ph.D. degree in Industrial and Manufacturing Systems Engineering from the University of Hong Kong, Hong Kong, in 2008.

He is a Full Professor with the School of Intelligent Systems Science and Engineering, Jinan University (Zhuhai Campus), Zhuhai, China. He has undertaken over 20 research projects funded by government and industry and has published nearly 200 technical papers, half of which have appeared in reputable journals. His research interests include IoT-based smart manufacturing systems, logistics and supply chain management, and industrial product/production service systems. 

Dr. Qu serves as the director or board member for several academic associations in the fields of industrial engineering and smart manufacturing.
\end{IEEEbiography}

\begin{IEEEbiography}[{\includegraphics[width=1in,height=1.25in,clip,keepaspectratio]{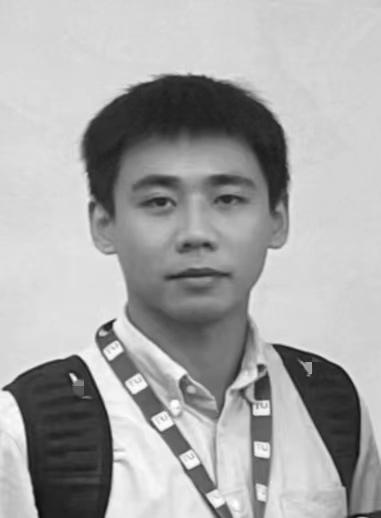}}]{Hainan Huang}  received his B.Sc. and M.Phil. degrees from the University of South China, Hengyang, China, in 2019 and 2022, respectively. He is currently pursuing a Ph.D. degree in Management Science and Engineering at Jinan University, Guangzhou, China. His current research interests include intelligent supply chain operations planning, decentralized supply chain planning, and distributed optimization.
\end{IEEEbiography}
\begin{IEEEbiography}
[{\includegraphics[width=1in,height=1.25in,clip,keepaspectratio]{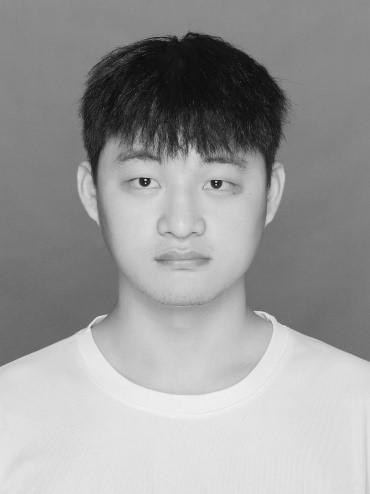}}]{Yafeng Wei}  received his B.Eng. degree in Information Management and Information Systems from Xinxiang Medical University, Xinxiang, China, in 2020, followed by an M.Sc. degree in Management Science and Engineering from Nanchang Institute of Technology, Nanchang, China, in 2023. He is currently pursuing a Ph.D. degree in Management Science and Engineering at Jinan University, Guangzhou, China. His current research interests include distributed optimization, intelligent manufacturing, and production scheduling.
\end{IEEEbiography}

\end{document}